\DeclareSymbolFont{AMSb}{U}{msb}{m}{n}
\documentclass[reqno,centertags,11pt,a4paper,noamsfonts]{amsart}
\pdfoutput=1
\usepackage[svgnames]{xcolor}
\usepackage{graphicx}
\usepackage[english]{babel}
\usepackage[babel=true,expansion=alltext,protrusion=alltext-nott,final]{microtype}
\usepackage[margin={3cm},marginparwidth={2.5cm}]{geometry}
\usepackage[utf8]{inputenc}
\usepackage{textcomp}
\usepackage[font=small]{subcaption}
\usepackage{algorithm}
\usepackage[algo2e]{algorithm2e}
\usepackage[sb]{libertine}
\usepackage[libertine,bigdelims,vvarbb]{newtxmath}
\useosf
\usepackage[varqu,varl]{inconsolata}
\usepackage[supstfm=libertinesups,%
       supscaled=1.2,%
       raised=-.13em]{superiors}
\usepackage{mathtools}
\usepackage[T1]{fontenc}
\usepackage{textcomp}
\usepackage{amsthm}
\usepackage[inline]{enumitem}
\numberwithin{equation}{section}
\usepackage[normalem]{ulem}

\usepackage{float}  
\usepackage{caption}
\usepackage{subcaption} 
\usepackage{xspace}         
\usepackage[scientific-notation=true]{siunitx}      
\usepackage{pstricks}
\usepackage{tikz}
\usetikzlibrary{positioning} 
\usetikzlibrary{calc}
\usepackage{pgfplots}
\pgfplotsset{width=10cm,compat=1.9}
\usepackage{bm}
\usepackage{sidecap}
\usepackage{graphicx,color}

\usepackage{amsmath}
\DeclareFontFamily{U}{mathx}{}
\DeclareFontShape{U}{mathx}{m}{n}{<-> mathx10}{}
\DeclareSymbolFont{mathx}{U}{mathx}{m}{n}
\DeclareMathAccent{\widehat}{0}{mathx}{"70}
\DeclareMathAccent{\widecheck}{0}{mathx}{"71}

\providecommand{\mr}[1]{\href{http://www.ams.org/mathscinet-getitem?mr=#1}{MR~#1}}
\providecommand{\zbl}[1]{\href{https://zbmath.org/?q=an:#1}{Zbl~#1}}

\providecommand{\doi}[2]{\href{https://doi.org/#1}{#2}}

\newcommand{\C}{\mathcal{C}}

\newcommand{\ii}{\imath}

\definecolor{light_gray}{gray}{0.75}
\definecolor{lighter_gray}{gray}{0.5}
\colorlet{light_blue}{blue!20}
\definecolor{dark_green}{rgb}{0.0, 0.6, 0.0}
\definecolor{royal_blue}{rgb}{0.0, 0.22, 0.66}
\definecolor{salmon}{rgb}{1.0, 0.55, 0.41}
\definecolor{gold}{rgb}{0.8, 0.63, 0.21}
\definecolor{navy_blue}{rgb}{0.0, 0.0, 0.5}
\definecolor{crimson}{rgb}{0.79, 0.0, 0.09}
\definecolor{amethyst}{rgb}{0.6, 0.4, 0.8}
\definecolor{alizarin}{rgb}{0.82, 0.1, 0.26}
\definecolor{amaranth}{rgb}{0.9, 0.17, 0.31}
\definecolor{azure}{rgb}{0.0, 0.5, 1.0}
\definecolor{canaryyellow}{rgb}{0.82, 0.41, 0.12}
\definecolor{carrotorange}{rgb}{0.8, 0.33, 0.0}
\definecolor{cadmiumgreen}{rgb}{0.0, 0.42, 0.24}
\definecolor{copper}{rgb}{0.72, 0.45, 0.2}
\definecolor{aqua}{rgb}{0.5, 1.0, 0.83}
\definecolor{awesome}{rgb}{1.0, 0.13, 0.32}
\definecolor{candyapplered}{rgb}{1.0, 0.03, 0.0}
\definecolor{caribbeangreen}{rgb}{0.0, 0.8, 0.6}
\definecolor{indigo}{rgb}{0.0, 0.25, 0.42}

\DeclareMathOperator{\weaklystar}{\rightharpoonup\kern-2.2ex ^* \, \,}

\def\XXint#1#2#3{{\setbox0=\hbox{$#1{#2#3}{\int}$ }
\vcenter{\hbox{$#2#3$ }}\kern-.6\wd0}}

\newcommand{\R}{\mathbb R}
\newcommand{\N}{\mathbb N}

\renewcommand{\C}{\mathbb C}

\newcommand{\bv}{{\bf v}}
\newcommand{\bu}{{\bf u}}

\newcommand{\bx}{{\bf x}}
\newcommand{\by}{{\bf y}}

\newcommand\norm[1]{\left\lVert #1 \right\rVert}
\newcommand\snorm[1]{\lVert #1 \rVert}

\newcommand\inner[1]{\langle #1 \rangle}

\newcommand{\ra}{\rightarrow}

\newcommand{\be}{\boldsymbol{e}}
\newcommand{\bof}{\boldsymbol{f}}
\newcommand{\bg}{\boldsymbol{g}}
\newcommand{\mP}{\mathbb{P}}
\newcommand{\mL}{\mathrm{L}}

\renewcommand{\phi}{\varphi}

\newcommand{\ext}{\mathrm{ext}}

\newcommand{\mH}{\mathrm{H}}
\newcommand{\mD}{\mathrm{D}}
\newcommand{\bold}[1]{{\bf #1}}
\newcommand{\ee}{\mathrm{e}}

\theoremstyle{plain}
\newtheorem{theorem}{Theorem}[section]
\newtheorem{proposition}[theorem]{Proposition}

\newtheorem{lemma}[theorem]{Lemma}

\newtheorem*{theorem*}{Theorem}

\theoremstyle{definition}
\newtheorem{notation}[theorem]{Notation}
\newtheorem{definition}[theorem]{Definition}
\newtheorem{remark}[theorem]{Remark}
\newtheorem*{remark*}{Remark}

\usepackage[pagebackref=false]{hyperref}
\hypersetup{
	linktoc=all,
	bookmarksdepth=3,
	colorlinks=true,
	linkcolor=Green,citecolor=Orange,urlcolor=DarkBlue,
	pdftitle={Trace estimates for harmonic functions along circular arcs and applications},
	pdfauthor={Thiago Carvalho Corso, Muhammad Hassan, Abhinav Jha, Benjamin Stamm}, 
	pdfsubject={},
	pdfkeywords={}} 
\begin{document}
\numberwithin{table}{section}
\title[Trace estimates for harmonic functions along circular arcs and applications]{Trace estimates for harmonic functions along circular arcs with applications to domain decomposition on overlapping disks}

\author[T.~Carvalho~Corso]{Thiago Carvalho Corso}
\address[T.~Carvalho Corso]{Institute of Applied Analysis and Numerical Simulation, University of Stuttgart, Pfaffenwaldring 57, 70569 Stuttgart, Germany}
\email{thiago.carvalho-corso@mathematik.uni-stuttgart.de}

\author[M.~Hassan]{Muhammad Hassan}
\address[M.~Hassan]{Institute of Mathematics, EPF Lausanne}
\email{muhammad.hassan@epfl.ch}

\author[A.~Jha]{Abhinav Jha}
\address[A.~Jha]{Department of Mathematics, Indian Institute of Technology Gandhinagar, Palaj, Gandhinagar, Gujarat-382355, India}
\email{abhinav.jha@iitgn.ac.in}

\author[B.~Stamm]{Benjamin Stamm}
\address[B.~Stamm]{Institute of Applied Analysis and Numerical Simulation, University of Stuttgart, Pfaffenwaldring 57, 70569 Stuttgart, Germany}
\email{Benjamin.Stamm@mathematik.uni-stuttgart.de}
\date{\today}
\thanks{\emph{Funding information}: DFG -- Project-ID 440641818; DFG -- Project-ID 442047500 -- SFB 1481.\\[1ex]
\textcopyright 2024 by the authors. Faithful reproduction of this article, in its entirety, by any means is permitted for noncommercial purposes.}

\begin{abstract}In this paper we derive several (and in many cases sharp) estimates for the $\mL^2$-trace norm of harmonic functions along circular arcs. More precisely, we obtain geometry-dependent estimates on the norm, spectral radius, and numerical range of the Dirchlet-to-Dirichlet (DtD) operator sending data on the boundary of the disk to the restriction of its harmonic extension along circular arcs inside the disk. The estimates we derive here have applications in the convergence analysis of the Schwarz domain decomposition method for overlapping disks in two dimensions. In particular, they allow us to establish a rigorous convergence proof for the discrete parallel Schwarz method applied to the Conductor-like Screening Model (COSMO) from theoretical chemistry in the two-disk case, and to derive error estimates with respect to the discretization parameter, the number of Schwarz iterations, and the geometry of the domain. Our analysis addresses challenges beyond classical domain decomposition theory, especially the weak enforcement of boundary conditions.
\end{abstract}
\maketitle
\tableofcontents
\setcounter{secnumdepth}{3}
\addtocontents{toc}{\protect\setcounter{tocdepth}{1}}

\section{Introduction}

This article is concerned with the study of the harmonic Dirichlet-to-Dirichlet (DtD) operator that takes as input $\mL^2$ boundary data on the two-dimensional disk, constructs the harmonic extension inside the disk, and outputs the trace along a circular arc inside the disk. More precisely, we explicitly characterize the spectrum and eigenfunctions of the composition of such DtD operators, and provide novel (and in many cases sharp) estimates on their operator norm and numerical range with respect to the $\mL^2$ trace inner-product. The derivation of these bounds is strongly motivated by applications in the convergence analysis of Schwarz domain decomposition methods, which have been applied in so-called implicit solvation models, and which we now briefly introduce.

\subsection{Motivation} Many phenomena in chemical physics involve the study of solute molecules embedded in a solvent, and it has been observed that certain properties of the solute molecule can change drastically with the choice of solvent. Attempts to correctly model solute-solvent interactions have led to the development of an extensive literature on so-called solvation models. Explicit solvation models consider an explicit description of both the solute molecule and the surrounding solvent molecules, thus leading to an accurate but computationally expensive description of solute-solvent systems. Implicit solvation models by comparison involve only solute degrees of freedom and describe the solvent as a dielectric continuum thereby resulting in a dramatic reduction in the computational complexity of the model. Within the class of implicit solvation models, the COnductor-like Screening MOdel (COSMO) \cite{KS93, TS95, BC98} is a popular choice due to its relative simplicity. Indeed, the COSMO equations consist of a simple Laplace boundary value problem of the form:
\begin{alignat}{3}\label{eq:poisson_problem}
    -\Delta u&=0 &&\mathrm{in}\ \ \Omega,\nonumber\\
    u&=g\ &&\mathrm{on}\ \ \partial\Omega,
\end{alignat}
where $g\in\mL^2(\partial\Omega)$ is some given Dirichlet boundary condition, $\Omega\subset \mathbb{R}^3$ is the solute cavity, and $\partial\Omega$ is the boundary of $\Omega$. Typically, the solute cavity is described as a union of $N$ open balls, i.e.,
$$
\Omega = \bigcup_{i=1}^N\Omega_i,
$$
with each $\Omega_i$ representing a spherical atomic cavity having boundary $\partial\Omega_i$.

The previous decade has seen significant development of domain decomposition methods for implicit solvation models such as the COSMO equations (see, e.e.,\cite{CMS13, SCLM16, QSM19}). In particular, the contribution \cite{CMS13} proposed the use of a parallel Schwarz-based Trefftz method for the solution of the COSMO equations \eqref{eq:poisson_problem}. This is a natural choice since the solute cavity $\Omega$ naturally partitions into an overlapping union of open balls, on which we have an explicit description of harmonic functions (e.g., in terms of Fourier series in 2D and spherical harmonics series in 3D). Introducing a discretisation parameter $L \in \mathbb{N}$ that limits the number of maximal terms in these infinite series, one arrives at a practical numerical method (see, e.g., Algorithm~\ref{alg:ddCOSMO} for a 2D version of this method).  

The resulting algorithm-- termed ddCOSMO-- has been successfully employed to obtain computationally efficient solutions to the COSMO implicit solvation model (see, e.g., \cite{CJL+15, LSL+14}). Independent of its practical utility, the ddCOSMO algorithm also raises several interesting questions from the analysis perspective, foremost among them whether it is possible to derive a \emph{rigorous proof of convergence} of Algorithm \ref{alg:ddCOSMO} as well as error estimates in terms of the discretization parameter $L$, the number of Schwarz iterations $k$ and the problem geometry $\Omega$. It is important to note that such convergence and error analyses do not immediately follow from the classical domain decomposition theory developed by Lions \cite{Lio88, Lio89} and are by no means obvious. 

\begin{algorithm}
\caption{Parallel Schwarz-based Trefftz Method for the Dirichlet Problem}\label{alg:ddCOSMO}
\hspace{-3cm}\KwData{Discretisation parameter $L$, boundary data $g \in \mL^2(\partial \Omega)$}
\hspace{-0.4cm}\KwResult{Approximate solution $u$ to the Dirichlet problem for boundary data $g$}
\hspace{-1.5cm}\textbf{Initialisation:} Select $u^{(0)} \in \mathcal{C}(\overline{\Omega})$ and set index $k=0$.\\
\While{convergence criterion not met}{
    $k = k + 1$\\
    \For{$j = 1$ to $N$}{
    \begin{itemize}
        \item Construct local boundary data $g_j\in \mL^2(\partial \Omega_j)$ by setting
        \[
        g_j = \begin{cases}
        g & \text{on } \partial \Omega_j \cap \partial \Omega\\
        u^{(k-1)} & \text{on } \partial \Omega_j \setminus \partial \Omega
        \end{cases}
        \]
       \item Find the Fourier coefficients $[g_j]_{{\ell}}$ for all $\vert \ell\vert \leq L$ given by
        \[
        [g_j]_{\ell} = \frac{1}{2\pi} \int_0^{2\pi} g_j(r_j e^{i \theta} + \mathbf{x}_j) e^{-i \ell \theta} \, d\theta,
        \]
        where the evaluation of $g_j$ is done by identifying $\C \simeq \R^2$.
        
        \item Define the harmonic function $u_j \in \mL^2(\Omega_j)$:
        \[
        u_j(\mathbf{x}) = [g_j]_0 + \sum_{\ell=1}^{L} \left( [g_j]_{\ell} (\mathbf{x} - \mathbf{x}_j)^{\ell} + [g_j]_{-\ell} \overline{(\mathbf{x} - \mathbf{x}_j)}^{\ell} \right),
        \]
        where $(\bx-\bx_j)^\ell$ is understood as the multiplication of complex numbers $\bx \in \R^2 \simeq \C$, and the bar denotes complex conjugation.
    \end{itemize}
    }
    Construct $u^{(k)} \in \mathcal{C}(\overline{\Omega})$ using local solutions $\{u_j\}_{j=1}^N$ and a partition of unity
}
\end{algorithm}

To see this, notice that Lions' framework of alternating projections \cite{Lio88} was developed for solutions to the Poisson equation with zero Dirichlet boundary conditions, i.e., for functions in the Sobolev space $\mH^1_0(\Omega)$, and the interpretation of the Schwarz method as a sequence of alternating projections crucially relies on the fact that for any $\Omega_i \subset \Omega$, we have the inclusion $\mH^1_0(\Omega_i) \subset \mH^1_0(\Omega)$. On the other hand, the governing equation in the present situation is the Laplace equation with $\mL^2$ Dirichlet boundary conditions so that the relevant function space is different. While this issue can be resolved by imposing additional regularity on the boundary data and transforming the Laplace equation to a homogeneous Poisson equation through a judicious lifting, a more fundamental obstacle is that the local sub-problems arising in Algorithm \ref{alg:ddCOSMO} are solved using a Trefftz method in which the PDE is satisfied exactly while the boundary conditions are imposed only in a weak sense. Consequently, even if one uses liftings, the relevant functions do not fit the alternating projection framework of Lions. Moreover, the fact that the local problems in Algorithm \ref{alg:ddCOSMO} are solved by imposing the boundary conditions only in a weak sense also prevents the use of the $\mL^{\infty}$-maximum principle arguments for harmonic functions for the convergence analysis (see Lions' second contribution on Schwarz method~\cite{Lio89}).

The existing articles on the analysis of ddCOSMO-type algorithms such as Algorithm \ref{alg:ddCOSMO} have therefore focused on the rate of convergence of the Schwarz iterates as a function of the problem geometry in the special case when the discretisation parameter $L=\infty$, i.e., when the local subproblems arising in Algorithm \ref{alg:ddCOSMO} are solved exactly. Note that in this situation, the boundary conditions appearing in the local sub-problems are enforced exactly, unlocking both $\mL^{\infty}$-maximum principle arguments and variational techniques for the convergence analysis. Contributions in this direction include the fundamental work of Ciaramella and Gander \cite{CG17, CG18a, CG18b}, which contains a convergence analysis of the parallel Schwarz method using the $\mL^{\infty}$-maximum principle for various geometries including the overlapping union of disks and shows that, in many cases, the convergence rate does not depend on the \emph{number of subdomains} in the geometry. Further work by Ciaramella, Hassan, and Stamm-- also based on $\mL^{\infty}$ maximum principle arguments \cite{CHS20a, CHS20b}-- studies the degradation in the convergence of the parallel Schwarz method for geometries consisting of the overlapping union of open balls with the additional property that some open balls lie entirely within the union of other balls (i.e., away from the boundary of the global domain). More recently, Reusken and Stamm \cite{RS21} have used variational arguments to perform a more careful convergence analysis of the parallel Schwarz method for generic geometries consisting of the union of overlapping balls and have obtained explicit dependencies of the convergence rate depending only on one parameter that reflects the globularity of the geometry.

\subsection{Main contribution} 
The primary goal of this article is to develop the tools necessary for the convergence analysis of the discrete parallel Schwarz-based Trefftz method for the Laplace Dirichlet problem when the set $\Omega$ is given by the union of overlapping disks. To this end, we shall derive several novel estimates on the $\mL^2$-trace norm along circular arcs of harmonic functions. More precisely, the main contribution of this paper can be summarized as follows:
\begin{itemize}
\item We present a detailed spectral analysis of the Dirchlet-to-Dirichlet (DtD) operator sending $\mL^2$ boundary data to the restriction of its harmonic extension along a circular arc inside the disk. In particular, we explicitly characterize the spectrum and eigenfunctions of the composition of such DtD operators, and provide novel (and in many cases sharp) estimates on their operator norm and numerical range with respect to the $\mL^2$ inner-product induced by the one-dimensional Hausdorff measure.
\end{itemize}
As an application of these estimates, we obtain the following results on the convergence of the ddCOSMO algorithm in the case of two disks:
\begin{enumerate}[label=(\roman*)]
\item We establish the well-posedness of the continuous (infinite-dimensional) ddCOSMO equation assuming only $\mL^2$ regularity on the boundary data. Moreover, in this case, we compute the sharp, geometry-dependent rate of convergence of the Schwarz iterates. 
\item We establish the well-posedness of the discrete (finite $L$) ddCOSMO equation. In this case, we also provide geometry-dependent, but possibly non-optimal, upper bounds on the asymptotic rate of convergence of the Schwarz iterates. 
\item We establish geometry-dependent error estimates between the solution of the discrete and continuous ddCOSMO problems in terms of the discretization parameter $L\in \N$.
\end{enumerate}

To our knowledge, these are the first rigorous results on the existence, convergence, and error analysis of the Schwarz-based Trefftz method for the Laplace equation on non-rectangular domains in the discrete case. A natural question, which is left open here, is whether the current analysis can be extended to the case of multiple overlapping disks in 3D, which is most relevant in applications. While we believe that the tools developed here can be extended to study the cases of multiple overlapping disks in 2D, the extension is not trivial and we postpone this development to future contribution. The extension to 3D, on the other hand, seems to require different tools, as most of our methods are based on complex analysis techniques.

\section{Main results}
\label{sec:mainresults} In this section, we present our main results concerning $\mL^2$ trace estimates for harmonic functions along circular arcs. We then show how these results can be used to obtain a convergence analysis for the ddCOSMO algorithm. 

\subsection{Trace estimates for harmonic functions along circular arcs}

To state our main estimates precisely, let us introduce the following notation. 

We denote by $\Omega \subset \R^2$ the open set
\begin{align*}
    \Omega := B_{r_1}(\bx_1) \cup B_{r_2}(\bx_2)
\end{align*}
where $B_{r_j}(\bx_j)$ is the open disk of radius $r_j$ and center $\bx_j$ in $\R^2$. We will assume that the intersection $B_{r_1}(\bx_1) \cap B_{r_2}(\bx_2)$ is non-empty, and that neither of the disks is fully contained on the other. In particular, the following sets are non-empty:
\begin{notation}[Interior and exterior boundaries]\label{not:boundaries} Let $\Omega= \Omega_1 \cup \Omega_2$ where $\Omega_i = B_{r_i}(\bx_i) \subset \R^2$ as defined above. Then for each $i \in \{1, 2\}$ we define
\begin{align*}
    \Gamma_{i}^{\rm ext} := \partial \Omega_i \cap \partial \Omega \quad \text{and} \quad \Gamma_{i}^{\rm int} := \partial \Omega_i  \setminus \Gamma_{i}^{\rm ext}.
\end{align*}
We refer to $\Gamma_{i}^{\rm ext} $ and $\Gamma_{i}^{\rm int}$ as the exterior and interior boundaries of $\Omega_i$ respectively.    
\end{notation}
\begin{definition}[Continuous local and global spaces on the boundary]\label{def:continuous spaces} We denote by $\mL^2(\partial \Omega_j)$ the usual space of square-integrable functions on the boundary with respect to the one-dimensional Hausdorff measure, i.e., 
\begin{align}\nonumber
    \mL^2(\partial \Omega_j) = \{ f: \partial \Omega_j \rightarrow \C &\colon \norm{f}_{\mL^2(\partial \Omega_j)}^2 = \inner{f,f}_{\mL^2(\partial \Omega_j)} <\infty \}, \qquad \text{where}\\[0.25em]
    \inner{f,g}_{\mL^2(\partial \Omega_j)} &= \int_{\partial \Omega_j} f(r) \overline{g(r)} \mathrm{d} \mathscr{H}^1(r), \label{eq:innerproduct}
\end{align}
and $\mathscr{H}^1$ denotes the one-dimensional Hausdorff measure. Moreover, we denote by $\mathcal{H}$ the Hilbert space given by the product (or direct sum)
\begin{align}
    \mathcal{H} \coloneqq \mL^2(\partial \Omega_1) \times  \mL^2(\partial \Omega_2) = \left\lbrace \bof = (f_1,f_2) : f_j \in \mL^2(\partial \Omega_j), j \in \{1,2\}  \right\rbrace,
\end{align}
endowed with the canonical inner product
\begin{align*}
    \inner{\bof, \bg}_{\mathcal{H}} = \inner{f_1,g_1}_{\mL^2(\partial \Omega_1)} + \inner{f_2, g_2}_{\mL^2(\partial \Omega_2)},
\end{align*}
where $\inner{f_j,g_j}_{\mL^2(\partial \Omega_j)}$ is defined according to \eqref{eq:innerproduct}. 
\end{definition}

We now introduce the Dirichlet-to-Dirichlet (DtD) map that plays a central role in this paper.
\begin{definition}[Harmonic Dirichlet-to-Dirichlet map]\label{def:trace} Let $\Omega$ be defined as above, and let the curve $\Gamma_j^{\rm int} = \partial \Omega_j \cap \Omega_{i}, ~ i, j \in  \{1, 2\}$ with $i \neq j$ be the interior part of the boundary of $\Omega_j$. Then, we define the \emph{Dirichlet-to-Dirichlet} ({\rm DtD}) map
\begin{align*}
    \gamma_{j} \colon \mL^2(\partial \Omega_{i}) \rightarrow \mL^2(\Gamma_j^{\rm int})
\end{align*}
as the mapping with the property that for any $g \in  \mL^2(\partial \Omega_{i})$, the output $\gamma_j g \in \mL^2(\Gamma_j^{\rm int}) $ is the pointwise restriction along $\Gamma_j^{\rm int}$ of the harmonic extension of $g$ on $\Omega_i$ given by the Poisson integral, i.e., $\gamma_j g = \mathrm{P}_i[g]\rvert_{\Gamma_j^{\rm int}}$, where 
\begin{align}
    \mathrm{P}_i[g](\bx) := \frac{1}{2\pi} \int_{\partial \Omega_i} \frac{r_i^2 - |\bx-\bx_i|^2}{r_i|\bx-\by|^2} g(\by) \mathrm{d} \mathscr{H}^1(\by) \qquad \mbox{for all }\bx \in \Omega_{i}, \label{eq:Poissonintegralball}
\end{align}
with $r_{i}$ and $\bx_i$ denoting, respectively, the radius and the center of the disk $\Omega_i$. 
\end{definition}

\begin{remark} [Extension by zero] 
The fact that $\gamma_j g$ belongs to $\mL^2(\Gamma_j^{\rm int})$ is a consequence of Theorem~\ref{thm:DtDests} below. Moreover, we shall often identify (without further notice) $\gamma_j g$ with the function on $\mL^2(\partial \Omega_j)$ given by extending $\gamma_j g$ by zero on the exterior part $\Gamma_j^{\rm ext}$.
\end{remark}

\begin{remark}[General framework for other PDEs] The Dirichlet-to-Dirichlet operator above encodes all the information on the partial differential equation (PDE) we wish to solve, namely, the Laplace equation in $\Omega$. In particular, the framework presented here could be easily extended to more general PDEs by simply replacing the harmonic Dirichlet-to-Dirichlet map with the one associated with the PDE of interest. Of course, for other PDEs, more suitable choices for the approximation space $\mathcal{H}_{L}$ (see Definition~\ref{def:continuous spaces} below) should be considered and the techniques required for the analysis of the corresponding DtD operator might be completely different. \end{remark}

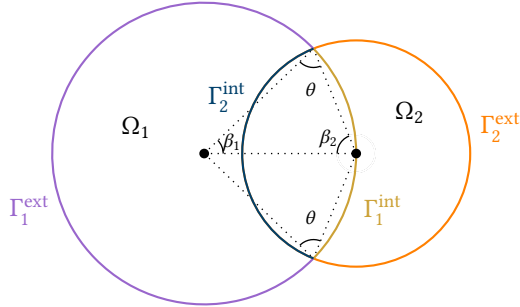
\begin{figure}[h!]
\begin{center}
\begin{tikzpicture}[scale=1.0]

\draw[line width=0.8pt, orange] (1,0) circle [radius=1.5];
\draw[line width=0.5pt] (1,0) circle [radius=0.25];
\draw[line width=0.8pt, white] (0.75, 0) arc(180:472.0243128:0.25);
\draw[line width=0.5pt, dotted] (1,0)--(0.4375	,1.3905372163) ;\draw[line width=0.5pt, dotted] (-1,0)--(1,0) ;
\draw[line width=0.5pt, dotted] (-1,0)--(0.4375	,1.3905372163); 
\draw[line width=0.5pt, dotted] (1,0)--(0.4375	,-1.3905372163) ;
\draw[line width=0.5pt, dotted] (-1,0)--(0.4375	,-1.3905372163); 
\draw[line width=0.8pt,amethyst] (0.4375	,1.3905372163) arc(44.0486257:315.9513743:2);
\draw[line width=0.8pt,gold] (0.4375	,-1.3905372163) arc(315.9513743:404.0486257:2);
\draw[line width=0.8pt,indigo] (0.4375	,1.3905372163) arc(112.0243128:247.9756872:1.5);
\filldraw[line width=0.8pt] (-1,0) circle [radius=0.05];
\filldraw[line width=0.8pt] (1,0) circle [radius=0.05];
\draw[line width=0.5pt] (0.258, 1.217) arc(224.1103974:291.8390154:0.25);
\draw[line width=0.5pt] (0.5313, -1.1588) arc(68.1609846:135.889603:0.25);
\draw[line width=0.5pt] (-0.75, 0) arc(0:44.0486257:0.25);
\begin{scriptsize}
\draw (-3.3,-0.75) node[align=center] {\small{$\textcolor{amethyst}{\Gamma^{\ext}_1}$}};
\draw (-1.9,0.35) node[align=center] {\small{$\Omega_1$}};
\draw (1.7,0.55) node[align=center] {\small{$\Omega_2$}};
\draw (-0.7,0.75) node[align=center] {\small{$\textcolor{indigo}{\Gamma_2^{\rm int}}$}};
\draw (2.9,0.35) node[align=center] {\small{$\textcolor{orange}{\Gamma_2^{\rm ext}}$}};
\draw (1.35,-0.75) node[align=center] {\small{$\textcolor{gold}{\Gamma^{\rm int}_1}$}};
\draw (0.4,0.85) node[align=center] {\tiny{$\theta$}};
\draw (0.4,-0.85) node[align=center] {\tiny{$\theta$}};
\draw (0.635,0.2) node[align=center] {\tiny{$\beta_2$}};
\draw (-0.635,0.15) node[align=center] {\tiny{$\beta_1$}};
\end{scriptsize}
\end{tikzpicture}
\caption{Visual example of two overlapping disks and associated angles.}\label{fig:main_def}
\end{center}
\end{figure}

Having introduced the DtD map, we now state the main technical contribution of this paper, whose proof can be found in Sections~\ref{sec:bipolartransform}--\ref{sec:proof}. 
\begin{theorem}[Estimates on the harmonic Dirichlet-to-Dirichlet map]\label{thm:DtDests} Let $\Omega_1$ and $\Omega_2$ be two disks as described above. Let $\theta \in (0,\pi)$ be the angle between the lines connecting the center of the disks to one of the intersection points of their boundaries, let $j \in \{1, 2\}$, and let $\beta_j$ be half the aperture angle of the interior arc $\Gamma^{\rm int}_j$, as depicted in Figure~\ref{fig:main_def}. Then the following holds:
\begin{enumerate}[label=(\roman*)]
\item (Operator norm) \label{it:operatornorm} The operator norm of $\gamma_j : \mL^2(\partial \Omega_i) \rightarrow \mL^2(\partial \Omega_j)$ satisfies the bound \begin{align}
    1 \leq \norm{\gamma_j} \leq  \sqrt{2}. \label{eq:DtDnorm}
\end{align}
\item\label{it:spectralradius} (Spectral radius) The spectral radius of the operator $\gamma_i \gamma_j : \mL^2(\partial \Omega_i) \rightarrow \mL^2(\partial \Omega_i)$ is given by
\begin{align}
    \label{eq:spectralradiusest} \rho(\gamma_i \gamma_j ) \coloneqq \sup \{|\lambda|: \lambda \in \sigma(\gamma_i \gamma_j) \} = \frac{1-\cos(\theta)}{2}.
\end{align}
Moreover, every $0 < \lambda < \frac{1-\cos(\theta)}{2}$ is an eigenvalue of $\gamma_i \gamma_j$.
\item\label{it:ext and  int} (Interior and exterior estimates) There exists $h^{\rm int} \in \mL^2(\Gamma_i^{\rm int})$ such that for any $g \in \mL^2(\Gamma_i^{\rm int})$ extended by zero to the rest of $\partial \Omega_i$, it holds 
\begin{align}
    \norm{\gamma_j g}_{\mL^2(\partial \Omega_j)}^2 \leq \sin(\theta/2)^2 \norm{P_{h^{\rm int}}^\perp g}_{L^2(\Gamma_i^{\rm int})}^2 + |\inner{h^{\rm int}, g}|^2, \label{eq:interiorbound}
\end{align}
where $P_{h^{\rm int}}^\perp$ is the $\mL^2$-orthogonal projection on the complement of $h^{\rm int}$. Similarly, there exists $h^{\rm ext} \in \mL^2(\Gamma_i^{\rm ext})$ such that for any $g \in \mL^2(\Gamma^{\rm ext}_i)$ extended by zero to the rest of $\partial \Omega_i$, it holds that
\begin{align}
    \norm{\gamma_j g}_{\mL^2(\partial \Omega_j)}^2 \leq \cos(\theta/2)^2 \norm{P_{h^{\rm ext}}^\perp g}_{\mL^2(\Gamma_i^{\rm ext})}^2 + |\inner{h^{\rm ext}, g}|^2. \label{eq:exteriorbound}
\end{align}

\item (Numerical radius) \label{it:numericalrange} We have
\begin{align}
    \frac{1+\sin(\theta/2)}{2} \leq \sup_{\substack{u_1 \in \mL^2(\partial \Omega_1)\\ u_2 \in \mL^2(\partial\Omega_2)}} \left\{ \frac{|\inner{u_1, \gamma_1 u_2}_{\mL^2(\partial \Omega_1)}+\inner{u_2, \gamma_2 u_1}_{\mL^2(\partial \Omega_2)}|}{\norm{u_1}_{\mL^2(\partial \Omega_1)}^2 + \norm{u_2}_{\mL^2(\partial \Omega_2)}^2} \right\} \leq f(\theta),\label{eq:numericalrange}
\end{align}
where $f(\theta)$ is positive and strictly smaller than one for all $\theta \in (0, \pi)$, and satisfies $f(\theta) = \frac{1+\sin(\theta/2)}{2}$ for $\theta$ large.

\end{enumerate}
\end{theorem}

\begin{remark}[Improved estimates] \label{rem:sharpness} Let us now make several remarks on the above estimates.
\begin{enumerate}[label=(\roman*)]
\item (Operator norm) The upper bound in~\eqref{eq:DtDnorm} is \emph{not} sharp and a better (but still not sharp) upper bound is given in Section~\ref{sec:proof} (see Theorem~\ref{thm:generalbound}). On the other hand, the lower bound is shown to be sharp as long as $\beta_j \leq \pi/2$, i.e., we have $\norm{\gamma_j} =1$ in this case. However, we shall show in Theorem~\ref{thm:generalbound} that for $\theta$ small and $\beta_j$ big, it holds that $\norm{\gamma_j} >1$. This result is rather surprising as it shows that the error after a single iteration of the Schwarz domain decomposition algorithm may in fact grow.
\item (Spectrum) We can, in fact, fully characterize the spectrum of $\gamma_i \gamma_j$ and explicitly compute its eigenfunctions (see Theorem~\ref{thm:spectrum}).  
\item (Interior and exterior estimates) \label{it:int ext remark} One can show that the operator norm of $\gamma_j$ when restricted to the interior, respectively, exterior parts of $\partial \Omega_i$ is always bounded from below by the constants $\sin(\theta/2)$ and $\cos(\theta/2)$ appearing in estimate~\eqref{eq:interiorbound} and~\eqref{eq:exteriorbound}. Moreover, for large values of $\theta$, one can show (see Theorem~\ref{thm:interior_exterior}) that $h^{\rm int} = 0$ and $h^{\rm ext}=0$, and therefore, these estimates are also sharp for certain configurations.  
\item (Numerical radius) The function $f$ in~\eqref{eq:numericalrange} is explicitly given by
\begin{align*}
    f(\theta) = \frac{1}{2} \left(1+\sin(\theta/2)^2 + g(\theta) + \sqrt{4\sin(\theta/2)^2 + h(\theta)}\right)^{\frac{1}{2}},
\end{align*}
where
\begin{align*}
\qquad \qquad &g(\theta) \coloneqq \left(\frac{\alpha(\theta)_+}{2} \frac{\theta}{\pi} +  \frac{\cos(\theta)_+}{2}\frac{\pi-\theta}{\pi}\right)\frac{\pi-\theta}{\pi} -\sin(\theta/2)^2\left(1 - \sqrt{1 + \frac{\alpha(\theta)_+}{\sin(\theta/2)^2} \frac{\theta}{\pi} \frac{\pi-\theta}{\pi}}\right), \\
        &\alpha(\theta) \coloneqq \frac{ \pi^2 \sin(\theta/2)^2 \cos(\theta) + \theta^2}{\pi^2 \sin(\theta/2)^2 + \theta^2 \cos(\theta)},  \\
	&h(\theta) \coloneqq g(\theta)^2 +  2g(\theta) \left(1+\sin(\theta/2)^2\right) - \cos(\theta/2)^2 \cos(\theta)_+ \frac{(\pi-\theta)^2}{\pi^2} , 
\end{align*}
and $a_+ = \max \{0,a\}$. A visual comparison between $f(\theta)$ and the lower bound in~\eqref{eq:numericalrange} can be seen in Figure~\ref{fig:f_function}. In particular, we note that $f(\theta) = \frac{1+\sin(\theta/2)}{2}$ for large values of $\theta$, and estimate~\eqref{eq:numericalrange} is therefore sharp on this range.
\end{enumerate}
\end{remark}

\input{plot_upper_lower_bound}

\subsection{Applications}
\label{sec:applications}
While we postpone its proof, we show in this section how Theorem~\ref{thm:DtDests} can be used to establish a complete error and convergence analysis for the ddCOSMO equations in the case of two disks. 

\subsubsection{The ddCosmo equations for two disks} 
We begin by introducing the notation and definitions necessary to state the ddCOSMO equations in the discrete (i.e., $L \in \mathbb{N}$) and continuous ($L=\infty$) setting. 

The approximate local and global function spaces used in Algorithm~\ref{alg:ddCOSMO} are the following.
\begin{definition}[Discrete local and global function spaces on the boundary]\label{def:discrete_approx} Let $\Omega_j = B_{r_j}(\bx_j)$ and $\mL^2(\partial \Omega_j)$ and $\mathcal{H}$ be the continuous spaces introduced in Definition~\ref{def:continuous spaces}. Then, for $L \in \N$, we denote by $\mathcal{H}_{L}(\partial \Omega_j) \subset {\mL}^{2}(\partial\Omega_j)$  the finite-dimensional subspace given by the span of the functions
\begin{align*}
    \left\{ (\bx-\bx_j)^\ell, \overline{(\bx-\bx_j)}^\ell\right\}_{\ell \leq L},
\end{align*}
where $\bx_j$ denotes the center of $\Omega_j$. In polar coordinates with respect to the center, these are simply the Fourier basis $\{r_j^{|\ell|} \ee^{\ii \ell \theta}\}_{|\ell| \leq L}$.
We also denote by $\mP_{j,L} : \mL^2(\partial \Omega_j) \rightarrow \mathcal{H}_{L}(\partial \Omega_j)$ the associated $\mL^2$-orthogonal projection on $\mathcal{H}_{L}(\partial \Omega_j)$. Moreover, we denote by $\mathcal{H}_{L} \subset \mathcal{H}$ the product space
\begin{align}
    \mathcal{H}_{L} = \mathcal{H}_{L}(\partial \Omega_1) \times \mathcal{H}_{L}(\partial \Omega_2),
\end{align}
 and by $\mP_{L}\colon \mathcal{H} \rightarrow \mathcal{H}_{L}$ the associated $\mathcal{H}$-orthogonal projection on $\mathcal{H}_{L}$.
\end{definition}

In other words, each $\mathcal{H}_L(\partial \Omega_j)$ is the space spanned by the restrictions to the boundary $\partial \Omega_j$ of all harmonic polynomials of degree up to $L$.

Using the DtD maps defined in Definition~\ref{def:trace}, we can introduce the so-called ddCOSMO operator.
\begin{definition}[The ddCOSMO operator] \label{def:operator}
Let $\Omega$ be as above, let $\mathcal{H}$ be the global space from Definition~\ref{def:continuous spaces}, and let $\gamma_j$, $j\in \{1,2\}$ be the Dirichlet-to-Dirichlet map introduced in Definition \ref{def:trace}. Then we define the operator $A\colon \mathcal{H} \rightarrow \mathcal{H}$ as
\begin{align}
    \bu = (u_1,u_2) \mapsto A\bu = (u_1-\gamma_1 u_2, u_2-\gamma_2 u_1) \in \mathcal{H}.\label{eq:Adef}
\end{align}
\end{definition}

We can now use the operator $A$ to state the continuous (infinite-dimensional) ddCOSMO equations, which Algorithm~\ref{alg:ddCOSMO} attempts to solve approximately, as well as the discrete ddCOSMO equations, which is the system of linear equations that the block-Jacobi methodology of Algorithm \ref{alg:ddCOSMO} is being applied to.
\vspace{2mm}

\noindent \textbf{Continuous ddCOSMO equations}~

\noindent Let  $\bg:= (g_1, g_2) \in \mathcal{H}$ and let $A$ be the operator given in Definition~\ref{def:operator}. We seek $\bu=(u_1, u_2) \in \mathcal{H}$ that solves the equation
\begin{align}\label{eq:continuous_ddCOSMO}
    A \bu = \bg.
\end{align}

\vspace{2mm}
\noindent \textbf{Discrete ddCOSMO equations}~ 

\noindent Let $\bg:= (g_1, g_2) \in \mathcal{H}$, let $L \in \mathbb{N}$, and let the approximation space $\mathcal{H}_{L} \subset \mathcal{H}$ and orthogonal projection operator $\mP_{L}$ be defined as in Definition \ref{def:discrete_approx}. We seek $\bu_{L}=(u_1^{L}, u_2^{L}) \in \mathcal{H}_{L}$ that solves the finite-dimensional equation
\begin{align}\label{eq:discrete_ddCOSMO}
    A_{L} \bu_{L} =  \mP_{L} \bg \qquad \text{where } ~A_{L} \coloneqq \mP_{L} A \mP_{L}.
\end{align}

\vspace{2mm}

Let us now briefly comment on the challenges in the analysis of the ddCOSMO equations. First, it is clear that the discrete ddCOSMO equation is simply a Galerkin discretization of the continuous ddCOSMO equation with discretization space $\mathcal{H}_{L} \subset \mathcal{H}$. Moreover, note that the operator $A$ can be written as
\begin{align}
    A = I- B,\label{eq:Adecomposition} 
\end{align}
where $I$ denotes the identity in $\mathcal{H}$ and $B : \mathcal{H} \rightarrow \mathcal{H}$ is defined as
\begin{align}
    B(u_1,u_2) \coloneqq (\gamma_1 u_2, \gamma_2 u_1). \label{eq:Bdef}
\end{align}
Thus, understanding the convergence of the ddCOSMO algorithm reduces to understanding the spectral properties of the operator $B$ and its reduced version $P_L B P_L$. This task, however, is non-trivial for several reasons. 

First, the operator $B$ is \emph{non-normal} and we can not immediately use tools from the spectral theory of normal operators. Second, one can show that the operator norm of $B$ is always greater or equal to $1$ (see Lemma~\ref{lem:Bests}). As $B$ is precisely the iteration operator of the Jacobi iterative procedure to solve the continuous ddCOSMO equations, this implies, in particular, that the error (measured in the $\mL^2$ norm) after a single iteration step can actually increase (see Remark~\ref{rem:sharpness}). Consequently, our analysis must not only include reasonable estimates on the operator norm of $B$ but also on its spectral radius, because the latter is directly related to the well-posedness of the problem and dominates the asymptotic convergence rate of the iterative procedure in the continuous setting. 

The main difficulty in the analysis of the discrete problem comes from the fact that the spectral radius of the discretized operator $\mathbb{P}_L B \mathbb{P}_L$ is not necessarily controlled by the spectral radius of $B$.  In particular, even with good estimates on the spectral radius of $B$, we can not immediately transfer the results from the continuous problem to the discrete one. A standard approach to overcoming this difficulty in the analysis of Galerkin discretizations is to establish the coercivity of the continuous operator, which is straightforwardly transferable to the discretized one. This is precisely the approach we follow here. Indeed, we observe that the coercivity constant of the operator $A$ is given by $1-\alpha(B)$, where $\alpha(B)$ denotes the spectral abscissa of $B$, i.e.,
\begin{align*}
    \alpha(B) \coloneqq \sup_{\bu \in \mathcal{H}} \frac{\mathrm{Re} \, \inner{\bu, B \bu}_{\mathcal{H}}}{\norm{\bu}_{\mathcal{H}}^2}.
\end{align*}Therefore, we can apply Theorem~\ref{thm:DtDests} to obtain lower bounds on the coercivity constant of $A$, and establish a full convergence and error analysis for both the discrete and the continuous ddCOSMO equations for two disks.

To carry out the above developments, the key result we shall use is the following corollary of Theorem~\ref{thm:DtDests}.
\begin{lemma}[Spectral properties of $B$]\label{lem:Bests} Let $\gamma_i$, $i \in \{1,2\}$ be the Dirichlet-to-Dirichlet map introduced in Definition~\ref{def:trace} and let $B: \mathcal{H} \rightarrow \mathcal{H}$ be the linear operator defined in~\eqref{eq:Bdef}. Then $B$ satisfies
\begin{enumerate}[label=(\roman*)]
\item (Operator norm)\label{it:operatornormB} The operator norm of $B$ satisfies \begin{align}\label{eq:Bnorm}
1 \leq \norm{B}_{\mathcal{H} \rightarrow \mathcal{H}} \leq \sqrt{2}
\end{align}
\item (Spectral radius)\label{it:spectralradiuesB} The spectral radius of $B$ is given by
\begin{align}
    \rho(B) = \sup \{ |\lambda| : \lambda \in \sigma(B)\} = \sin(\theta/2), \label{eq:spectralradiuesB}
\end{align}
where $\theta$ is the angle between the lines connecting the center of each disk to one of the intersection points of their boundaries.
\item (Numerical radius) \label{it:coercivity} $B$ satisfies the estimate
\begin{align}|\inner{B \bu,\bu}| \leq f(\theta) \norm{\bu}_{\mathcal{H}}^2, \quad \mbox{for any $\bu \in \mathcal{H}$,} \label{eq:coercivity}
\end{align}
where $f(\theta)$ is the function given in Remark~\ref{rem:sharpness}.
\end{enumerate}
\end{lemma}

\begin{proof}
    For the boundedness estimate, we note that by estimate~\eqref{eq:DtDnorm}, we have
    \begin{align*}
        \norm{B \bu}_{\mathcal{H}}^2 &= \sum_{j\neq i} \norm{\gamma_i u_j}_{\mL^2(\partial \Omega_i)}^2 \leq 2\sum_{j=1}^2 \norm{u_j}_{\mL^2(\partial \Omega_j)}^2 = 2\norm{\bu}_{\mathcal{H}}^2, \quad \mbox{for any $\bu \in \mathcal{H}$,}
    \end{align*}
    which proves the upper bound $\norm{B}^2 \leq 2$. The lower bound follows from the lower bound in~\eqref{eq:DtDnorm}. The estimate on the numerical radius of $B$ is also immediate from the definition of $B$ and estimate~\eqref{eq:numericalrange} in Theorem~\ref{thm:DtDests}.

To compute the spectral radius of $B$ we note that $B^2$ can be written as
\begin{align*}
    B^2 = \begin{pmatrix} \gamma_1 \gamma_2 & 0 \\ 0 & \gamma_2\gamma_1\end{pmatrix}
\end{align*}
in block representation with respect to the decomposition $\mathcal{H} = \mL^2(\partial \Omega_1) \oplus \mL^2(\partial \Omega_2)$. Hence $\sigma(B^2) = \sigma(\gamma_1 \gamma_2) \cup \sigma(\gamma_2 \gamma_1)$ and from~\eqref{eq:spectralradiusest} we have
\begin{align*}
    \rho(B^2) = \frac{1-\cos(\theta)}{2} = \sin(\theta/2)^2.
\end{align*}
To conclude, note that by Gelfand's formula we have 
\begin{align*}
    \rho(B) = \lim_{n \ra \infty} \norm{B^n}^{\frac{1}{n}} = \left(\lim_{n \ra \infty} \norm{B^{2n}}^{\frac{1}{n}}\right)^{\frac{1}{2}} = \rho(B^2)^{\frac{1}{2}} = \sin(\theta/2).
\end{align*}
 \end{proof}

\subsubsection{Well-posedness of the ddCOSMO equations} We now prove the unique solvability of the continuous problem \eqref{eq:continuous_ddCOSMO} and the Galerkin discretization \eqref{eq:discrete_ddCOSMO}. 

\begin{theorem}[Well-posedness of the ddCOSMO equations] \label{thm:wellposedness} For any $\bg\in \mathcal{H}$ there exists a unique solution $\bu \in \mathcal{H}$ of the continuous ddCOSMO equation~\eqref{eq:continuous_ddCOSMO}. Moreover, for any $L \in \N$, there exists a unique solution $\bu^L \in \mathcal{H}_L$ satisfying the discrete ddCOSMO equations~\eqref{eq:discrete_ddCOSMO}.
\end{theorem}

\begin{proof} As $\rho(B) = \sin(\theta/2) <1$ by Lemma~\ref{lem:Bests}~\ref{it:spectralradiuesB}, the operator $A = I-B$ is invertible and the continuous problem is well-posed. For the discrete problem, we note that from Lemma~\ref{lem:Bests}~\ref{it:coercivity} and the Cauchy-Schwarz inequality, 
\begin{align}
    \norm{A_L \bu}\norm{\mathbb{P}_L \bu} \geq |\inner{A_L \bu , \bu}| = |\norm{\mathbb{P}_L \bu}^2 + \inner{B \mathbb{P}_L\bu \mathbb{P}_L \bu}| \geq \left(1- f(\theta)\right) \norm{\mathbb{P}_L \bu}^2 \label{eq:coercivityofA0},
\end{align}
where $f(\theta) <1$. In particular, $A_L$ is a bounded-below mapping on $\mathcal{H}_L$, and hence injective, and therefore, invertible.
\end{proof}

\begin{remark}[Optimal bound on coercivity] \label{rem:coercivityA} Estimate~\eqref{eq:coercivityofA0} gives a lower bound on the coercivity constant of $A$ in terms of the function $f(\theta)$ introduced in Remark \ref{rem:sharpness}. 
Moreover, we note that, since the DtD operator has a positive integral kernel (the Poisson kernel), the spectra abscissa and numerical radius of $B$ coincide. Consequently, the lower bound in~\eqref{eq:numericalrange} yields an upper bound on the coercivity of $A$, which together with the lower bound in~\eqref{eq:coercivityofA0} reads
\begin{align*}
    1-f(\theta) \leq \inf_{u \in \mathcal{H}\setminus \{0\} } \frac{\mathrm{Re} \inner{u, Au}}{\norm{u}^2}  \leq \frac{1-\sin(\theta/2)}{2} \label{eq:coercivityofA}.
\end{align*}
In particular, estimate~\eqref{eq:numericalrange} yields the sharp value for the coercivity constant of $A$ for $\theta$ close to $\pi$.
 \end{remark}
\subsubsection{Convergence rate of the parallel Schwarz iterates} We now provide upper bounds on the asymptotic rate of convergence for the parallel Schwarz domain decomposition method in terms of the number of iterations for both the continuous ($L =\infty$) and discrete ($L \in \mathbb{N}$) case. To the best of our knowledge, this is the first rigorous proof of convergence for the ddCOSMO algorithm in the discrete setting.  In the continuous setting, previous results exist under the assumption that the boundary date is either bounded or has higher regularity in the Sobolev scale (see, e.g., \cite{CG18a, Lio88, RS21}). However, a proof of convergence with explicit upper bounds on the convergence rate under the sole assumption that the boundary data is in $\mL^2(\partial \Omega)$ is also new.

Recalling the definitions of $A$ and $B$, we see that algorithm~\ref{alg:ddCOSMO} is simply the Jacobi iterative method applied to the operator $A$ (respectively, $\mP_{L} A \mP_{L}$ in the discrete case) for solving the linear equations $A \bu = \bg$. In particular, at the $n^{th}$ iteration, the approximate solution $\bu^{(n)}$ is given by
\begin{align*}
    \bu^{(n)}_{L} = \mP_{L} \bg + \mP_{L} B \bu^{(n-1)}_{L}. 
\end{align*}
Hence the error term $\be^{(n)}_{L} \coloneqq  \bu_{L} - \bu^{(n)}_{L}$ satisfies the equation
\begin{align*}
    \be^{(n)}_{L} = \mP_{L} B  \be^{(n-1)}_{L}.
\end{align*}
If each local sub-problem is solved exactly, i.e., $L = \infty$, the projection $\mathbb{P}_L$ becomes the identity. In particular, we can obtain the \emph{optimal} convergence rate of the continuous problem by using Lemma~\ref{lem:Bests}.

\begin{proposition}[Optimal convergence rate of parallel Schwarz iterates in the continuous case]\label{lem:convergencecontinuous} Let the operator $B$ be defined as in~\eqref{eq:Bdef} and let $\be^{(n)} = B^n \be^{(0)}$ for some $\be^{(0)} \in \mathcal{H}$. Then for any $\epsilon>0$ we have
\begin{align}
    \snorm{\be^{(n)}} \leq \left(\sin(\theta/2) +\epsilon\right)^n \snorm{\be^{(0)}}
\end{align}
for any $n$ large enough. In particular $\bu^{(n)} = \bg + B \bu^{(n-1)}$ converges exponentially fast to the solution $\bu$ of the continuous ddCOSMO equation $A \bu = \bg$. Moreover, the above rate is sharp in the sense that for any $\epsilon>0$, there exists $\be^{(0)} \in \mathcal{H}$ such that
\begin{align}
    \snorm{\be^{(n)}} > \left(\sin(\theta/2)-\epsilon\right)^{n} \snorm{\be^{(0)}}, \quad \mbox{for all $n \in \N$.} \label{eq:lowerboundconvergencerate}
\end{align}
\end{proposition}
\begin{proof} The proof of the upper bound is immediate from Gelfand's formula for the spectral radius
\begin{align*}
    \rho(B) = \lim_{n \ra \infty} \left( \snorm{B^n}\right)^{\frac{1}{n}},
\end{align*}
and Lemma~\ref{lem:Bests}~\ref{it:spectralradiuesB}. The proof of the lower bound follows from the fact that every $0<\lambda < \sin(\theta/2)^2$ is an eigenvalue of $\gamma_1 \gamma_2$ (cf. Theorem~\ref{thm:DtDests}~\ref{it:spectralradius}). Indeed, let $f\in \mL^2(\partial \Omega_1)$ be an eigenfunction of $\gamma_1 \gamma_2$ with eigenvalue $\lambda = (\sin(\theta/2)-\epsilon)^2$, then $\be^{(0)} \coloneqq (f, \gamma_2 f/ \lambda^{\frac12}) \in \mathcal{H}$ is an eigenfunction of $B$ with eigenvalue $\lambda^{\frac12} = \sin(\theta/2) -\epsilon$ and the result follows.
\end{proof}

\begin{remark}[Non-asymptotic convergence rate] \label{rem:nonasymptotic} We emphasize that the lower bound in~\eqref{eq:lowerboundconvergencerate} is non-asymptotic, i.e., it holds for every $n \in \N$. Moreover, for $\theta$ large enough (e.g. $\theta \geq 3\pi/4$), it follows from Theorem~\ref{thm:interior_exterior} that the norm of $\gamma_1$ and $\gamma_2$ when restricted to the interior boundaries are given by $\sin(\theta/2)$. Consequently, $\norm{B^n} \leq \sin(\theta/2)^{n-1}$ and the following \emph{non-asymptotic} upper-bound holds:
\begin{align*}
    \snorm{\be^{(n)}} \leq \sin(\theta/2)^{n-1}\snorm{\be^{(0)}}, \quad \mbox{for any $n \in \N$.}
\end{align*}
\end{remark}

For the discrete problem we have the following upper bound on the asymptotic convergence rate of the Schwarz iterates. 

\begin{proposition}[Upper bound on convergence rate of the Schwarz iterates in the discrete case]\label{lem:convergencediscrete} Let the operator $B$ be defined as in~\eqref{eq:Bdef} and let $\be^{(n)} = (\mathbb{P}_L B \mathbb{P}_L)^n \be^{(0)}$ for some $\be^{(0)} \in \mathcal{H}$. Then for any $\epsilon>0$ and $n \in \N$ large enough, we have
\begin{align}
    \snorm{\be^{(n)}} \leq \left(f(\theta) +\epsilon\right)^n \snorm{\be^{(0)}},
\end{align}
where $f(\theta)$ is the function in Remark~\ref{rem:sharpness}. In particular, the sequence $\bu^{(n)} = \mathbb{P}_L \bg + \mathbb{P}_L B \mathbb{P}_L \bu^{(n-1)}$ converges exponentially fast to the solution $\bu \in \mathcal{H}_L$ of the discrete ddCOSMO equation $A_L \bu = \mathbb{P}_L \bg$. 
\end{proposition}

\begin{proof} As in the proof of Proposition~\ref{lem:convergencecontinuous}, it suffices to estimate the spectral radius of $B_L =\mathbb{P}_L B \mathbb{P}_L$. Hence, the result follows from the simple inequalities
\begin{align*}
\rho(B_L) \leq r(B_L) \leq r(B) \leq f(\theta),
\end{align*}
where $r(B)$ (respectively $r(B_L)$) denotes the numerical radius
\begin{align*}
    r(B) \coloneqq \sup_{0 \neq \bu \in \mathcal{H}} \frac{|\inner{\bu , B \bu}|}{\norm{\bu}^2},
\end{align*}
and the last inequality comes from Lemma~\ref{lem:Bests}~\ref{it:coercivity}.
\end{proof}

\begin{remark}[Degenerating convergence rate] Proposition~\ref{lem:convergencecontinuous} shows that, as the disks are pulled apart, the convergence rate decreases. This happens because the information flow between the two interior parts of the boundary increases, or more precisely, the composite operator $\gamma_1 \gamma_2$ converges strongly to the identity operator on the interior boundary, as the disks are pulled apart. Therefore, any error made in the initialization is propagated essentially unchanged at each iteration, which leads to a worse convergence rate. This is a well-known phenomena for overlapping domain decomposition methods. \end{remark} 

\subsubsection{Solution of the Laplace equation}

We have now proved that both the continuous and the discrete ddCOSMO equations are well-posed. The goal of this section is then to relate the unique solution of the continuous ddCOSMO equation with the solution of the underlying Laplace equation. More precisely, we seek to show that if $\bg:= (g_1, g_2) \in \mathcal{H}$ satisfies
\begin{align*}
    g_i \rvert_{\Gamma_i^{\rm int}} = 0, \quad i\in \{1,2\},
\end{align*}
then the unique solution $\bu \coloneqq A^{-1} \bg =(u_1, u_2) \in \mathcal{H}$ of the continuous ddCOSMO equation satisfies $u_1 = u\rvert_{\partial \Omega_1}$ and $u_2 = u\rvert_{\partial \Omega_2}$, where $u$ is the unique solution of the boundary value problem (BVP)
\begin{align}
\begin{dcases} \Delta u = 0, \quad \mbox{in $\Omega = \Omega_1 \cup \Omega_2$,} \\
u = g_1, \quad \mbox{on $\Gamma_1^{\rm ext}$,}\\
u = g_2, \quad \mbox{on $\Gamma_2^{\rm ext}$.} \end{dcases} \label{eq:LaplaceTwinDisk}
\end{align}

To do so, the typical approach in the domain decomposition literature is to show that the restriction of the solution to Equation \eqref{eq:LaplaceTwinDisk} belongs to the appropriate function space and satisfies the associated local equations (i.e., the continuous ddCOSMO equation in our case) whose solution is known to be unique. Under the assumption that the boundary data is either bounded or has higher Sobolev regularity, this is an immediate consequence of standard results in the well-posedness of problem~\eqref{eq:LaplaceTwinDisk}. However, under the sole assumption that the boundary data belongs to $\mL^2(\partial \Omega)$, we need to appeal to a more technical well-posedness result.

To state this well-posedness result precisely, let us first recall the definition of the Sobolev scale of spaces on a domain $\Omega$.
\begin{definition}[Sobolev spaces] \label{def:Sobolev} Let $\Omega \subset \R^2$ be a open subset. Then for any $s\geq 0$, we say that $u: \Omega \rightarrow \C$ belongs to $\mH^s(\Omega)$ if there exists a function $v: \R^2 \rightarrow \C$ satisfying $v\rvert_{\Omega} = u$ and
\begin{align*}
    \norm{v}_{\mH^s(\R^2)}^2 = \int_{\R^2} |\widehat{v}(\xi)|^2 (1+|\xi|)^s \mathrm{d} \xi < \infty,
\end{align*}
where $\widehat{v}$ denotes the Fourier transform of $v$. Moreover, we set
\begin{align*}
    \norm{u}_{\mH^s(\Omega)} = \inf_{v\rvert_{\Omega} = u} \norm{v}_{\mH^s(\R^2)}.
\end{align*}
\end{definition}
\begin{remark} For Lipschitz $\Omega$ and integer $s$, the above definition agrees with the standard definition of Sobolev spaces as the space of functions whose weak derivatives up to order $s$ belong to $\mL^2(\Omega)$ (see, e.g., \cite[Chapter 3]{McL00}).
\end{remark}

Using the above definition of Sobolev spaces, we can now state the low-regularity well-posedness result\footnote{The formulation of Theorem~\ref{thm:Dahlberg} in the Sobolev space $\mH^{\frac12}(\Omega)$ is equivalent to the original formulation of \cite{Dah77} in terms of the nontangential maximal function by the equivalence in \cite[Corollary 5.5]{JK95}.} for the BVP~\eqref{eq:LaplaceTwinDisk}.

\begin{theorem}[\cite{Dah77}] \label{thm:Dahlberg} Let $\Omega \subset \R^2$ be a Lipschitz domain. Then for any $g \in \mL^2(\partial \Omega)$ there exists a unique harmonic function $u \in \mH^{\frac12}(\Omega)$ such that $u\rvert_{\partial \Omega} = g$ in the nontangential limit sense, i.e., 
\begin{align*}
    \lim_{ y \ra \Gamma_{\partial \Omega}(\bx)} u(\by) = g(\bx) \quad \mbox{for $\mathscr{H}^1$-almost every $\bx\in \partial \Omega$},
\end{align*}
where $\Gamma_{\partial \Omega}(\bx)$ is the nontangential approach region\footnote{The factor of $2$ appearing in the definition of $\Gamma_{\partial \Omega}(\bx)$ is rather arbitrary; it could be replaced by any number greater than $1$.} with vertex $\bx \in \partial \Omega$, defined as
\begin{align*}
    \Gamma_{\partial \Omega}(\bx) \coloneqq \{ \by \in \Omega : |\bx-\by| < 2\, \mathrm{dist}(\by, \partial \Omega) \}.
\end{align*}
\end{theorem}
We can now apply Theorem~\ref{thm:Dahlberg} to prove the claim that the ddCOSMO solution is indeed the (restriction to $\partial \Omega_1 \cup \partial \Omega_2$ of the) solution of the BVP~\eqref{eq:LaplaceTwinDisk}.
\begin{theorem}[Solution of BVP via ddCOSMO] \label{thm:solution} Let $\bg = (g_1,g_2) \in \mathcal{H}$ satisfy
\begin{align*}
    g_i \rvert_{\Gamma_i^{\rm int}} = 0, \quad \mbox{for $i \in \{1,2\}$},
\end{align*}
and let $u \in \mH^{\frac12}(\Omega)$ be the unique solution of 
\begin{align*}
    \begin{dcases} -\Delta u = 0, \quad &\mbox{in $\Omega$,}\\
    u = g_1, \quad &\mbox{on $\Gamma_1^{\rm ext}$,} \\
    u = g_2, \quad &\mbox{on $\Gamma_2^{\rm ext}$,} \end{dcases}
\end{align*}
in the sense of Theorem~\ref{thm:Dahlberg}. Then the unique solution $\bu \coloneqq A^{-1} \bg$ of the ddCOSMO equation~\eqref{eq:continuous_ddCOSMO} satisfies
\begin{align*}
    u_1 = u \rvert_{\partial \Omega_1} \quad \mbox{and} \quad u_2 = u\rvert_{\partial \Omega_2}.
\end{align*}
\end{theorem}

\begin{proof} Since the domain $\Omega_1 \cap \Omega_2$ is Lipschitz, general trace estimates as in \cite[Theorem 3.1]{BGM22} guarantee that $u\rvert_{\Gamma_i^{\rm int}} \in \mL^2(\Gamma_i^{\rm int})$ for any $i \in \{1,2\}$. In particular, the functions
 \begin{align}
     v_i \coloneqq \begin{dcases} g_i, \quad &\mbox{on $\Gamma_i^{\rm ext}$,} \\ u \rvert_{\Gamma_i^{\rm int}}, \quad &\mbox{on $\Gamma_i^{\rm int}$.} \end{dcases} \label{eq:v_idef}
 \end{align}
 belong to $\mL^2(\partial \Omega_i)$ for any $i\in \{1,2\}$. Since $A$ is invertible (by Theorem~\ref{thm:wellposedness}), it is enough to show that $\bv = (v_1,v_2)$ satisfies $A \bv = \bg$ to complete the proof.
 
To this end, we first observe that, by Theorem~\ref{thm:Dahlberg}, the restriction $u \rvert_{\Omega_i}$ belongs to the Sobolev space $\mH^{\frac12}(\Omega_i)$ (see Definition~\ref{def:Sobolev}), and the nontagential limit of $u$ along $\Gamma_i^{\rm ext}$ is given by $g_i$. Moreover, as $u$ is smooth inside $\Omega$, we have
 \begin{align*}
     \lim_{\by \in \Gamma_{\partial \Omega_i}(\bx) \ra \bx} u(\by) = u(\bx), \quad \mbox{for any $\bx \in \Gamma_i^{\rm int} \setminus (\partial \Omega_1 \cap \partial \Omega_2)$,}
 \end{align*}
 where $\Gamma_{\partial \Omega_i}(\bx) = \{ \by \in \Omega_i : |\bx-\by| < 2\mathrm{dist}(\by,\partial \Omega_i)\}$. Thus $u\rvert_{\partial \Omega_i}$ is the unique solution of the BVP
 \begin{align}
     \begin{dcases} -\Delta u = 0, \quad \mbox{in $\Omega_i$,} \\
     u = v_i, \quad \mbox{on $\partial \Omega_i$} \end{dcases} \quad\mbox{with $v_i$ defined in \eqref{eq:v_idef},} \label{eq:BVP1}
 \end{align}
 in the sense of Theorem~\ref{thm:Dahlberg}.
 
 Next, observe that for any $g \in \mL^2(\partial \Omega_j)$, the Poisson integral $\mathrm{P}_j[g]$ defined according to \eqref{eq:Poissonintegralball} belongs to $\mH^{\frac12}(\Omega_j)$ and its nontangential limit along $\partial \Omega_j$ is given by $g$ (see \cite{Dah79}). In particular, by the uniqueness of the solution to \eqref{eq:BVP1}, we have
 \begin{align*}
     u\rvert_{\Omega_i} = \mathrm{P}_i[v_i] \quad \mbox{with $v_i$ defined in \eqref{eq:BVP1} for $i \in \{1,2\}$.}
 \end{align*}
 Therefore, by the definition of the Dirichlet-to-Dirichlet map we have
 \begin{align*}
     \gamma_i v_j = \mathrm{P}_j[v_j] \rvert_{\Gamma_i^{\rm int}}= u \rvert_{\Gamma_i^{\rm int}} = v_i \rvert_{\Gamma_i^{\rm int}}\quad \mbox{for $i \neq j \in \{1,2\}$.}
     \end{align*}
     We thus conclude that the function $\bv = (v_1,v_2) \in \mathcal{H}$ satisfies 
     \begin{align*}
         A \bv = (v_1 - \gamma_1 v_2, v_2 -\gamma_2 v_1) = (v_1 \rvert_{\Gamma^{\rm ext}_1}, v_2 \rvert_{\Gamma^{\rm ext}_2}) = (g_1,g_2) = \bg \in \mathcal{H},
     \end{align*}
     which completes the proof.
\end{proof}

\subsubsection{Error estimates} We can now combine the previous results with the classical trace theorem to obtain the following error estimates between the discrete and continuous ddCOSMO solutions in terms of the discretization parameter $L$.
\begin{theorem}[Discrete-to-continuum error estimates]\label{thm:errorestimates} Let $g \in \mL^2(\partial \Omega)$ and suppose that the solution of the BVP~\eqref{eq:LaplaceTwinDisk} satisfies $u \in \mH^{k+\frac12}(\Omega)$ for some $k \geq 0$. Let $\bu \in \mathcal{H}$ and $\bu_{L} \in \mathcal{H}_{L}$ be respectively the solutions of the continuous and discrete ddCOSMO equations with right hand side
\begin{align*}
    \bg \coloneqq (g\rvert_{\Gamma_{1}^{\rm ext}},g\rvert_{\Gamma_2^{\rm ext}}) \in \mathcal{H}.
\end{align*}
Then, we have
\begin{align}
    \norm{\bu-\bu_{L}}_{\mathcal{H}} \lesssim \frac{1}{1-f(\theta)} \frac{1}{L^k} \norm{u}_{\mH^{k+\frac12}(\Omega)},
\end{align}
where $f(\theta)$ is the function in Remark~\ref{rem:sharpness}, and $\theta$ is the angle between the lines connecting the center of each disk to one of the intersection point of their boundaries.
\end{theorem}

\begin{proof} As $A$ is coercive with coercivity constant $1-f(\theta)$ and is continuous with continuity constant $\norm{A} \leq 1 + \sqrt{2}$, from Cea's lemma we have
\begin{align*}
    \norm{\bu - \bu_L}_{\mathcal{H}} \leq \frac{1+\sqrt{2}}{1-f(\theta)} \norm{\mathbb{P}_L^\perp \bu}_\mathcal{H} = \frac{1+\sqrt{2}}{1-f(\theta)}\sqrt{\snorm{\mathbb{P}_{1,L}^\perp u_1}_{\mL^2(\partial \Omega_1)}^2 + \snorm{\mathbb{P}_{2,L}^\perp 
 u_2}_{\mL^2(\partial \Omega_2)}^2}.
\end{align*}
Since $u_1 = u\rvert_{\partial \Omega_1}$ and $u_2 = u\rvert_{\partial \Omega_2}$ by Theorem~\ref{thm:solution}, the result now follows from 
\begin{align*}
    \norm{\mathbb{P}_{L,j} u_j}_{\mL^2(\partial \Omega_j)} \lesssim \frac{1}{L^k} \norm{u_j}_{H^k(\partial \Omega_j)} \lesssim \frac{1}{L^k} \norm{u}_{H^{k+\frac12}(\Omega_j)},
\end{align*}
where the first inequality is well-known\footnote{And can easily be shown by using the Fourier series representation on the circle} and the second inequality follows from the classical trace theorem (see, e.g., \cite[Section 2.6]{Sau11}).
\end{proof}

We end this subsection by remarking that, from a computational point-of-view, one is often interested in the solution $u\in \mH^{\frac12}(\Omega)$ of the BVP~\eqref{eq:LaplaceTwinDisk}. However, once the solution $\mathbf{u}_L \in \mathcal{H}_L$ to the discrete ddCOSMO equations \eqref{eq:discrete_ddCOSMO} has been obtained, we can reconstruct the harmonic extension inside $\Omega_1$ and $\Omega_2$ exactly. Moreover, in view of Theorem \ref{thm:errorestimates}, these extensions will converge to $u$  (in $H^{\frac12}(\Omega_j)$) as the discretisation parameter $L \to \infty$ with the rates depending on the regularity of $u$ as stated in Theorem~\ref{thm:errorestimates}.

\subsection{Outline of the proof of Theorem~\ref{thm:DtDests}}\label{sec:mainproofoutline}

The remaining sections are devoted to the proof of Theorem~\ref{thm:DtDests}, which is the main result of this paper. As this proof is rather involved and interesting in its own right, we now outline the key steps and the main difficulties, and how these are distributed throughout the next sections.

The first step in the proof of Theorem~\ref{thm:DtDests} is to apply the bipolar transform, a conformal transformation between the punctured (with two holes) plane and the strip, to re-state the sought-after estimates for $\gamma_j$ in terms of estimates along horizontal lines for the Poisson integral on the strip. This will be done in Section~\ref{sec:bipolartransform}. The main advantage of working in the strip is that the DtD map restricted to each part (interior and exterior) becomes a convolution operator in $\R$ with an explicit kernel, namely, the Poisson kernel on the strip. Hence, the next natural idea is to apply the Fourier transform to turn the Poisson integral into a multiplication operator. However, this can not be done immediately as the underlying function spaces are weighed $\mL^2$ spaces with exponentially decaying weights and therefore contain functions that grow exponentially fast and have no well-defined Fourier transform. 

To overcome this difficulty, we use the symmetry of the Poisson kernel and a duality argument. More precisely, the second step of the proof, which is carried out in Section~\ref{sec:duality}, consists of using a duality argument to re-state the sought-after estimates on the weighed $\mL^2$ spaces as estimates on the dual weighted spaces. The dual spaces have exponentially growing weights, and therefore, only contain functions that decay exponentially fast. As a result, we can identify the Fourier transform of the dual weighted space with the $\mL^2$ Hardy space of holomorphic functions on a different (dual) strip in the complex plane. This allow us to look at the adjoint DtD map as a multiplication operator on the Hardy space.

The third and final step of the proof then consists in effectively estimating various spectral properties of this multiplication operator. This will be done in Section~\ref{sec:proof} and is the most involved (calculation heavy) part of the proof. The main tools used at this step are a combination of some standard complex and real analysis techniques, and some auxiliary lemmas, whose proofs are provided in Appendix~\ref{app:Poisson}.

\section{From the disk to the strip}\label{sec:bipolartransform}

In this section we show that the DtD operator can be seen as a convolution operator on certain weighted $\mL^2$ spaces along horizontal lines on the strip. This is the first step in the proof of Theorem~\ref{thm:DtDests}. 

\subsection{The bipolar transform} 

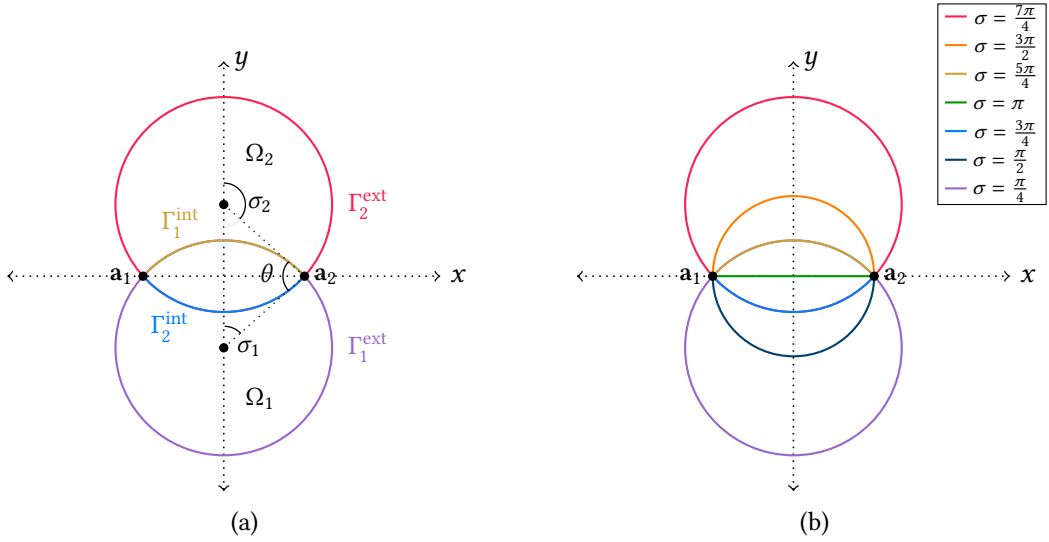
\begin{figure}[h!]
	\centering
	\begin{subfigure}{0.5\textwidth}
		\centering
		\begin{tikzpicture}[scale=0.95]
\draw[<->,line width=0.6pt, dotted] (-3,0)--(3,0) node[right]{$x$};
\draw[<->,line width=0.6pt, dotted] (0,-3)--(0,3) node[right]{$y$};
\draw[line width=0.8pt, color = amethyst] (0,-1) circle [radius=1.5];
\draw[line width=0.8pt, awesome] (0,1) circle [radius=1.5];

\draw[line width=0.8pt,gold] (1.1180339887,0) arc(41.8103149:138.1896851:1.5);
\draw[line width=0.8pt,azure] (-1.1180339887,0) arc(221.8103149:318.1896851:1.5);
\filldraw[line width=0.8pt] (1.1180339887,0) circle [radius=0.05];
\filldraw[line width=0.8pt] (-1.1180339887,0) circle [radius=0.05];
\filldraw[line width=0.8pt] (0,-1) circle [radius=0.05];
\filldraw[line width=0.8pt] (0,1) circle [radius=0.05];
\draw[line width=0.5pt, dotted] (0,-1)--(1.1180339887,0) ;
\draw[line width=0.5pt, dotted] (0,1)--(1.1180339887,0) ;
\draw[line width=0.5pt] (0, 0.7) arc(270:450:0.3);
\draw[line width=0.7pt, white] (0, 0.7) arc(270:318.3024534:0.3);
\draw[line width=0.5pt] (0.236, -0.8) arc(41.8122623:90:0.3);
\draw[line width=0.5pt] (0.897, 0.202) arc(137.6751427:222.3248573:0.3);
\begin{scriptsize}
\draw (2.0,-1) node[align=center] {\small{$\textcolor{amethyst}{\Gamma_1^{\ext}}$}};
\draw (0.5,1.7) node[align=center] {\small{$\Omega_2$}};
\draw (0.5,-1.7) node[align=center] {\small{$\Omega_1$}};
\draw (-0.75,-0.7) node[align=center] {\small{$\textcolor{azure}{\Gamma_2^{\rm int}}$}};
\draw (2.0,1) node[align=center] {\small{$\textcolor{awesome}{\Gamma_2^{\rm ext}}$}};
\draw (1.4180339887,0) node[align=center] {\small{$\mathrm{\bold{a}_2}$}};
\draw (-1.4180339887,0) node[align=center] {\small{$\mathrm{\bold{a}_1}$}};
\draw (-0.6,0.75) node[align=center] {\small{$\textcolor{gold}{\Gamma_1^{\rm int}}$}};
\draw (0.6,0.05) node[align=center] {\small{$\theta$}};
\draw (0.5,1) node[align=left] {\small{$\sigma_2$}};
\draw (0.35,-1) node[align=left] {\small{$\sigma_1$}};

\end{scriptsize}
\end{tikzpicture}
		\caption{ }		
		\label{fig:01a}
	\end{subfigure}\hfill
	\begin{subfigure}{0.5\textwidth}
		\centering
		\begin{tikzpicture}[scale=0.95]
\draw[<->,line width=0.6pt, dotted] (-3,0)--(3,0) node[right]{$x$};
\draw[<->,line width=0.6pt, dotted] (0,-3)--(0,3) node[right]{$y$};
\draw[line width=0.8pt, color = amethyst] (0,-1) circle [radius=1.5];
\draw[line width=0.8pt, awesome] (0,1) circle [radius=1.5];

\draw[line width=0.8pt,gold] (1.1180339887,0) arc(41.8103149:138.1896851:1.5);
\draw[line width=0.8pt,azure] (-1.1180339887,0) arc(221.8103149:318.1896851:1.5);
\draw[line width=0.8pt,orange] (1.1180339887,0) arc(0:180:1.1180339887);
\draw[line width=0.8pt,indigo] (-1.1180339887,0) arc(180:360:1.1180339887);
\draw[line width=0.8pt, dark_green] (-1.1180339887,0)--(1.1180339887,0);
\filldraw[line width=0.8pt] (1.1180339887,0) circle [radius=0.05];
\filldraw[line width=0.8pt] (-1.1180339887,0) circle [radius=0.05];

\draw[line width=0.8pt,amethyst] (2.1,1.225)--(2.4,1.225);
\draw[line width=0.8pt,indigo] (2.1,1.625)--(2.4,1.625);
\draw[line width=0.8pt,azure] (2.1,2.025)--(2.4,2.025);
\draw[line width=0.8pt,dark_green] (2.1,2.425)--(2.4,2.425);
\draw[line width=0.8pt,gold] (2.1,2.825)--(2.4,2.825);
\draw[line width=0.8pt,orange] (2.1,3.225)--(2.4,3.225);
\draw[line width=0.8pt,awesome] (2.1,3.625)--(2.4,3.625);
\draw (2,1)--(3.5,1)--(3.5,3.8)--(2,3.8)--cycle; 
\begin{scriptsize}
\draw (2.875,1.2) node[align=center] {$\sigma = \frac{\pi}{4}$};
\draw (2.87	5,1.6) node[align=center] {$\sigma = \frac{\pi}{2}$};
\draw (2.95,2) node[align=center] {$\sigma = \frac{3\pi}{4}$};
\draw (2.85,2.4) node[align=center] {$\sigma = \pi$};
\draw (2.95,2.8) node[align=center] {$\sigma = \frac{5\pi}{4}$};
\draw (2.95,3.2) node[align=center] {$\sigma = \frac{3\pi}{2}$};
\draw (2.95,3.6) node[align=center] {$\sigma = \frac{7\pi}{4}$};
\draw (1.4180339887,0) node[align=center] {\small{$\mathrm{\bold{a}_2}$}};
\draw (-1.4180339887,0) node[align=center] {\small{$\mathrm{\bold{a}_1}$}};
\end{scriptsize}
\end{tikzpicture}
		\caption{ }
		\label{fig:01b}
	\end{subfigure}
        \caption{(a) An example of the geometric setting and (b) curves of constant $\sigma$ in Cartesian space.}
	\label{fig:01}
\end{figure}

Let $\Omega_1$ and $\Omega_2$ be two open disks as described in Section~\ref{sec:mainresults}. Then the starting point of our analysis is the observation that all estimates in Theorem~\ref{thm:DtDests} are invariant under translations, rotations, and dilations of the chosen coordinate system. Therefore, without loss of generality, we can assume that the intersection points $\{\mathbf{a}_1, \mathbf{a}_2\} = \partial \Omega_1 \cap \partial \Omega_2$ are given by
\begin{align*}
	\mathbf{a}_1 := (-1, 0)\quad \mbox{and}\quad \quad \mathbf{a}_2 := (1, 0),
\end{align*}
and that the exterior part of the boundary of $\Omega_1$, denoted by $\Gamma_1^{\rm ext}$, is located on the lower half-plane $\{ (x,y) \in \R^2 : y\leq 0 \}$. An example of this geometric setting is displayed in Figure \ref{fig:01a}.

  Based on this observation, we introduce the so-called bipolar coordinate system with poles (or foci) at $\{\mathbf{a}_1, \mathbf{a}_2\}:=\{(-1,0),(1,0)\}$.

\begin{definition}[Bipolar Transform]\label{def:Bipolar} Let $\mathcal{S} = \R \times (0,2\pi) \subset \mathbb{R}^2$ denote the infinite horizontal strip. We define the bipolar transform as the mapping $\Psi \colon \mathcal{S} \rightarrow \mathbb{R}^2$ given by
	\begin{align}
		\Psi(\tau, \sigma) :=\left( \frac{\sinh(\tau)}{\cosh(\tau)-\cos(\sigma)}, \frac{-\sin(\sigma)}{\cosh(\tau)-\cos(\sigma)}\right) \qquad \forall (\tau, \sigma) \in \mathcal{S}. \label{eq:bipolardef}
	\end{align}
	\end{definition}

The bipolar transform has a number of useful properties that can, for instance, be found in \cite{MF53}. We list two in particular that will be crucial for our subsequent analysis. 

\begin{proposition}[The bipolar transform]\label{prop:bipolartransform} The bipolar transform $\Psi \colon \mathcal{S} \rightarrow \mathbb{R}^2$ defined through Definition \ref{def:Bipolar} has the following properties:
\begin{enumerate}[label=(\roman*)]
\item $\Psi$ is a conformal\footnote{In the sense that the Jacobian $\mD \Psi (\tau, \sigma)$ at any $(\tau, \sigma) \in \mathcal{S}=\R^2 \times (0, 2\pi)$ is the product of a non-zero (not necessarily positive) scalar and an orthogonal matrix.} diffeomorphism from the strip $\mathcal{S}$ to the set 
\begin{align*}
    \mathbb{R}^2 \setminus \left\{ (x,0) \in \R^2 : |x| \geq 1 \right\}.
\end{align*}
	
\item \label{it:disks} For any $\sigma\in (0, \pi)$, the bipolar transform $\Psi$ maps the open strip $\mathcal{S}_\sigma \coloneqq \R \times (\sigma,\sigma+\pi)$ to the unique open disk $\Omega \subset \R^2$ such that $\{\mathbf{a_1},\mathbf{a_2}\} \subset \partial \Omega$ and the center is located at
\begin{align*}
    \bx_\Omega = (0, -\cot\sigma).
\end{align*}
Moreover, the lines $\R_\sigma = \R \times \{\sigma\}$ and $\R_{\sigma+\pi}= \R \times \{\sigma+\pi\}$ are mapped to the portions of the boundary $\partial \Omega$ in the lower and upper half-plane, i.e.,
\begin{align*}
		\partial \Omega = \underbrace{\{(x, y) \in \partial \Omega \colon y < 0\}}_{=\Psi(\mathbb{R}_\sigma)} \cup \underbrace{\{(x, y) \in \partial \Omega \colon y >0\}}_{= \Psi(\mathbb{R}_{\sigma+\pi})} \cup \{\mathbf{a}_1, \mathbf{a}_2\},
	 \end{align*}
(see Figure~\ref{fig:01b}). In particular, for any open disk $\Omega \subset \R^2$ satisfying $\{\mathbf{a}_1, \mathbf{a}_2 \}$ $\subset \partial\Omega$ there exists $\sigma \in (0,\pi)$ such that $\Omega = \Psi(\mathcal{S}_\sigma)$. 
	\end{enumerate}
\end{proposition}

  In view of item~\ref{it:disks} of Proposition~\ref{prop:bipolartransform}, we introduce the following weighted Lebesgue spaces.

\begin{definition}[Weighted Lebesgue Spaces]\label{def:Leb_w} Let $\Psi \colon \mathcal{S} \rightarrow \mathbb{R}^2$ be the bipolar transform introduced in Definition~\ref{def:Bipolar}, let $\sigma \in (0,2\pi)$, let $p\in [1,\infty]$, and let $\Psi_{\sigma} \colon \R_\sigma \rightarrow \Psi(\R_\sigma)$ be the restriction $\Psi_\sigma = \Psi \rvert_{\R_{\sigma}}$ where $\R_\sigma := \R \times \{\sigma\}$. Then we define the weighted Lebesgue space $\mathcal{L}^p(\R_{\sigma})$ as the vector space 
	\begin{align*}
		\mathcal{L}^p(\R_{\sigma})	:=\{f \colon \mathbb{R}_{\sigma}\rightarrow \mathbb{C} ~\text{ such that } f\circ \Psi_{\sigma}^{-1}\in \mL^p(\Psi(\mathbb{R}_{\sigma}))\},
	\end{align*}
	endowed with the norm $\Vert f \Vert_{\mathcal{L}^p(\R_{\sigma})} :=\Vert f  \circ \Psi^{-1} \Vert_{\mL^p(\Psi(\mathbb{R}_{\sigma}))}$, where $\mL^p\bigr(\Psi(\R_\sigma)\bigr)$ denotes the $\mL^p$ space with respect to the one-dimensional Hausdorff measure on the circular arc $\Psi(\R_\sigma)$. Moreover, we denote by $\mathcal{L}^p_\sigma$ the direct sum space
 \begin{align*}
     \mathcal{L}^p_\sigma = \mathcal{L}^p(\R_\sigma) \oplus \mathcal{L}^p(\R_{\sigma+\pi}) = \mathcal{L}^p(\R_\sigma) \times \mathcal{L}^p(\R_{\sigma+\pi})
 \end{align*}
 endowed with the canonical norm
 \begin{align*}
     \norm{(f,g)}_{\mathcal{L}^p_\sigma}^p = \norm{f}_{\mathcal{L}^p(\R_\sigma)}^p + \norm{g}_{\mathcal{L}^p(\R_{\sigma+\pi})}^p.
 \end{align*}
\end{definition}

\begin{remark}[Alternative definition]\label{def:Leb_w_2} From the change of variables formula we can equivalently define $\mathcal{L}^p(\R_\sigma)$ as the vector space of (Lebesgue) measurable functions $f:\R \rightarrow \C$ satisfying
\begin{align}
     \Vert f \Vert_{\mathcal{L}^p(\R_{\sigma})}^p  :=  \int_{-\infty} ^{\infty} \vert f(\tau)\vert^p \frac{1}{\cosh(\tau) - \cos(\sigma)}\, d\tau < \infty .\label{eq:weightednorm}
\end{align}
In the sequel, we will frequently make use of this alternative representation of the weighted Lebesgue space $\mathcal{L}^p(\R_\sigma)$ without further comment.
\end{remark}

An immediate consequence of Proposition~\ref{prop:bipolartransform}~\ref{it:disks} is that the pullback map $\Psi^\#$, when restricted to the boundary of $\Omega_1$ (or $\Omega_2$), defines an isometry between the spaces $\mL^p(\partial \Omega_1)$ and $\mathcal{L}^p_\sigma$ for suitable $\sigma \in (0,\pi)$. More precisely,
we have the following statement.
\begin{lemma} \label{lem:bipolarisomorphism} Let $\Omega_1$ and $\Omega_2$ be two open disks with non-empty intersection and such that $\partial \Omega_1 \cap \partial \Omega_2 = \{\mathbf{a}_1, \mathbf{a}_2\}$ and $\Gamma_1^{\rm ext}$ lies on the lower half-plane. Then there exists $\sigma_1, \sigma_2 \in (0,\pi)$ such that $\sigma_1 < \sigma_2$ and
\begin{align*}
    \Gamma_1^{\rm ext} = \Psi(\R_{\sigma_1}), \quad \Gamma_1^{\rm int} = \Psi(\R_{\sigma_1+\pi}), \quad \Gamma_2^{\rm int} = \Psi(\R_{\sigma_2}), \quad\mbox{and}\quad \Gamma_2^{\rm ext} = \Psi(\R_{\sigma_2+\pi}).
\end{align*}
In particular, the pullback map $\Psi_j^\# : \mL^p(\partial \Omega_j) \rightarrow \mathcal{L}^p_{\sigma_j}$ defined as
\begin{align}
    g \in \mL^p(\partial \Omega_j)  \mapsto \Psi_j^\# g = \bigr(g\left(\Psi\left(\cdot,\sigma_j\right)\right), g\left(\Psi\left(\cdot, \sigma_j+\pi\right)\right)\bigr) \in \mathcal{L}^p(\R_{\sigma_j}) \times \mathcal{L}^p(\R_{\sigma_j+\pi}) \label{eq:pullback}
\end{align}is an isometric bijection for any $1\leq p \leq \infty$.
\end{lemma}
\begin{proof} The existence of $\sigma_1,\sigma_2 \in (0,\pi)$ follows from Proposition~\ref{prop:bipolartransform}~\ref{it:disks}. Since we assumed $\Gamma_1^{\rm ext} = \Psi(\R_\sigma)$ to be in the lower half-plane, we must have $\sigma_1 < \sigma_2$. Since the complement 
\begin{align*}
    \{\mathbf{a}_1 , \mathbf{a}_2\} = \Omega_j \setminus \left(\Gamma_j^{\rm int} \cup \Gamma_j^{\rm ext}\right)
\end{align*}
has (one-)Hausdorff measure zero, the proof follows from the canonical decomposition $\mL^p(\partial \Omega_j) = \mL^p(\Gamma_j^{\rm int}) \oplus \mL^p(\Gamma_j^{\rm int}) = \mL^p(\Gamma_j^{\rm int}) \times \mL^p(\Gamma_j^{\rm ext})$ and the definition of the weighted $\mL^p$ spaces.
\end{proof}

\begin{remark}[$\beta$ vs. $\sigma$ angles] Due to the convention adopted for the disks $\Omega_1$ and $\Omega_2$ and for the bipolar transform (see Figure~\ref{fig:01}), we have $\sigma_1 = \beta_1$ and $\sigma_2 = \pi-\beta_2$ where $\beta_j$ are the (half) aperture angles of the interior boundaries $\Gamma_j^{\rm int}$ depicted in Figure~\ref{fig:main_def}. \end{remark}
\subsection{The Poisson integral on the strip} We now observe that, since $\Psi$ is a conformal transformation, it maps harmonic functions on the disk $\Omega_j$ to harmonic functions on the strip $\mathcal{S}_{\sigma_j}$. Consequently, we can give an explicit description of the DtD harmonic map $\gamma_j$ in the weighted spaces $\mathcal{L}^2_\sigma$ by using a well-known Poisson representation formula for harmonic functions on the strip. To state this result precisely, let us introduce the Poisson kernel on the strip and recall the associated Poisson representation formula.
\begin{definition}[Poisson kernel on the infinite strip]\label{def:Poisson_strip} The Poisson kernel on the infinite strip $\mathcal{S}_0:= \R \times (0,\pi)$ is the function $\mathscr{P} \colon \mathcal{S}_0 \rightarrow \mathbb{R}$ given by
	\begin{align}
	\forall (\tau,\theta) \in \mathcal{S}_0\colon \quad	\mathscr{P}(\tau,\theta) = \mathscr{P}_\theta(\tau) := \frac{1}{2\pi}  \frac{\sin(\theta) }{\cosh(\tau)- \cos(\theta)}. \label{eq:Poissondef}
		\end{align}
	\end{definition}

\begin{lemma}[Poisson integral on the strip \cite{Wid61}]\label{lem:Widder} Let $\mathcal{S}_0:= \mathbb{R} \times (0, \pi)$ be the the infinite strip, let $\mathscr{P}\colon \mathcal{S}_0 \rightarrow \R$ be the Poisson kernel defined in \eqref{eq:Poissondef}. Then for any $g,h \in C_c^\infty(\R)$, the function
	\begin{align}
	 u_{g,h}(\tau,\theta) := \int_{-\infty}^{\infty} \mathscr{P}_\theta(\tau - \xi)g(\xi)\, d\xi+  \int_{-\infty}^{\infty} \mathscr{P}_{\pi- \theta}(\tau-\xi)h(\xi)\, d\xi \label{eq:poissonintegralstrip}
		\end{align}
is the unique bounded harmonic function in $\mathcal{S}_0$ satisfying $u_{g,h}(\cdot,\theta) \rightarrow g $ and $u_{g,h}(\cdot,\theta) \rightarrow h$ pointwise as $\theta \rightarrow 0$ and $\theta \ra \pi$, respectively. 
	\end{lemma}

\begin{proof}
The proof of the above lemma is rather standard and relies on the fact that the Poisson kernel is an approximation of the identity. For the details we refer, e.g., to \cite{Wid61}, where a detailed study of the Poisson kernel on the strip is carried out.
\end{proof}

The explicit description of $\gamma_j$ on the strip, which is the main result of this section, can now be stated as follows.
\begin{lemma}[Dirichlet-to-Dirichlet map on the strip]\label{lem:DtDonStrip} Let $\Omega_1$ and $\Omega_2$ be two open disks as described in Section~\ref{sec:mainresults} and $\gamma_j : \mathrm{L^2}(\partial \Omega_i) \rightarrow \mathrm{L^2}(\partial \Omega_j)$ be the associated DtD map introduced in Definition~\ref{def:trace}. Then there exists a conformal transformation $\Psi : \mathcal{S} \rightarrow \R^2 \setminus (\partial \Omega_1 \cap \partial \Omega_2)$ and $\sigma_1 < \sigma_2 \in (0,\pi)$ such that the pullback map $\Psi_j^\# : \mL^2(\partial \Omega_j) \rightarrow \mathcal{L}^2_{\sigma_j}$ (defined as in~\eqref{eq:pullback}) is an isometric bijection. Moreover, the composite maps
\begin{align*}
    \Psi_{1}^{\#} \circ \gamma_1 \circ (\Psi_{2}^\#)^{-1} : \mathcal{L}^2_{\sigma_2} \rightarrow \mathcal{L}^2_{\sigma_1} \quad \mbox{and}\quad \Psi_{2}^{\#} \circ \gamma_2 \circ (\Psi_{1}^\#)^{-1} : \mathcal{L}^2_{\sigma_1} \rightarrow \mathcal{L}^2_{\sigma_2}
\end{align*}
are given by
\begin{align}
    &(f, g) \mapsto \Psi_{1}^{\#} \circ \gamma_1 \circ (\Psi_{2}^\#)^{-1}(f,g) = \left(0,\,\,\,P_{\pi-\theta} f + P_{\theta} g\right), \label{eq:DtDstrip1} \\ 
    &(f, g) \mapsto \Psi_{2}^{\#} \circ \gamma_2 \circ (\Psi_{1}^\#)^{-1}(f,g) = \left(P_{\theta} f + P_{\pi-\theta} g,\,\,\, 0\right) \label{eq:DtDstrip2}
\end{align}
where $\theta \coloneqq \sigma_2-\sigma_1 \in (0,\pi)$ and $P_\theta$ is the convolution operator
\begin{align}
    (P_\theta g)(\tau) \coloneqq (\mathscr{P}_\theta \ast g)(\tau) = \frac{1}{2\pi} \int_\R \frac{\sin(\theta)}{\cosh(\xi) - \cos(\theta)} g(\tau-\xi) \mathrm{d} \xi. \label{eq:convolutionPoisson}
\end{align}
\end{lemma} 

\begin{proof} As previously mentioned, it suffices to consider the case where $\partial \Omega_1 \cap \partial \Omega_2 = \{ \mathbf{a}_1, \mathbf{a}_2\}$ and $\Gamma_1^{\rm ext}$ lies on the lower half-plane. In this case, we already know from Lemma~\ref{lem:bipolarisomorphism} that the Bipolar transform is a conformal transformation such that $\Psi^\#_{\sigma_j}: \mL^2(\partial \Omega_j) \rightarrow \mathcal{L}^2_{\sigma_j}$ is an isometric bijection. Hence, it remains to show that the composite maps are given by the formulas~\eqref{eq:DtDstrip1} and~\eqref{eq:DtDstrip2}. 

Since the proof of~\eqref{eq:DtDstrip1} is completely analogous, let us focus on~\eqref{eq:DtDstrip2}. In this case, we first note that from general trace estimates for harmonic functions (cf. \cite[Theorem 3.1]{BGM22}), the map $\gamma_2$ is bounded in  $\mL^2(\partial \Omega_1) \rightarrow \mL^2( \partial \Omega_2)$, and therefore, the composite map $\Psi_{\sigma_2}^{\#} \circ \gamma_2 \circ (\Psi_{\sigma_1}^\#)^{-1}$ is bounded from $\mathcal{L}^2_{\sigma_1}$ to $\mathcal{L}^2_{\sigma_2}$. Since $C_c^\infty(\R)$ is dense in $\mathcal{L}^2(\R_\sigma)$ for any $\sigma \in (0,2\pi)$\footnote{See the proof of Lemma~\ref{lem:residuethm} for a proof of this statement.}, it suffices to prove~\eqref{eq:DtDstrip1} in the case where $f = (f_{\sigma_1},f_{\sigma_1+\pi}) \in C_c^\infty(\R)\times C_c^\infty(\R)$. For such $f$, we note that by the conformal property of $\Psi$ and Lemma~\ref{lem:Widder}, the maps 
\begin{align*}
    (\tau,\sigma) \in \R \times(\sigma_1, \sigma_1+\pi) \mapsto u(\tau,\sigma) \coloneqq  P_{\Omega_1}[(\Psi_1^\#)^{-1} f]\left(\Psi(\tau,\sigma)\right) 
\end{align*}
and
\begin{align*}
    (\tau,\sigma) \in \R \times(\sigma_1, \sigma_1+\pi) \mapsto \widetilde{u}(\tau,\sigma) \coloneqq (P_{\sigma-\sigma_1} f_{\sigma_1})(\tau) + (P_{\pi-\sigma+\sigma_1} f_{\sigma_1+\pi})(\tau)
\end{align*}
are harmonic and bounded in $\mathcal{S}_{\sigma_1} = \R \times(\sigma_1,\sigma_1+\pi)$ and satisfy 
\begin{align*}
    \lim_{\sigma \uparrow \sigma_1+\pi} \!\! u(\tau,\sigma) = f_{\sigma_1+\pi}(\tau) =  \!\! \lim_{\sigma \uparrow \sigma_1+\pi} \widetilde{u}(\tau,\sigma) \quad\mbox{and}\quad \lim_{\sigma \downarrow \sigma_1} u(\tau,\sigma) = f_{\sigma_1}(\tau) = \lim_{\sigma \downarrow \sigma_1} \widetilde{u}(\tau,\sigma), \quad \mbox{for any $\tau \in \R$.}
\end{align*}
The result then follows from the uniqueness of such harmonic functions (Lemma~\ref{lem:Widder}) and the identity
\begin{align*}
    u(\tau,\sigma_2) = \gamma_2\left((\Psi_1^\#)^{-1}g\right)\left(\Psi(\tau,\sigma_2)\right) = \left(\Psi_2^\# \circ \gamma_2 \circ ( \Psi_1^\#)^{-1} f\right)(\tau).
\end{align*}
\end{proof}

\begin{remark}[Weighted $\mL^p$ estimates] \label{rem:Lpests} Lemmas~\ref{lem:bipolarisomorphism} and~\ref{lem:Widder} can also be used to formulate a $\mL^p$ version of Lemma~\ref{lem:DtDonStrip} in terms of the weighted Lebesgue spaces $\mathcal{L}^p_{\sigma}$ with $1\leq p \leq \infty$. In particular, the $\mL^\infty$ estimates obtained in \cite{CG18a} are equivalent to the operator norm estimate
\begin{align*}
    \norm{P_\theta}_{\mL^\infty(\R) \rightarrow \mL^\infty(\R)} = \frac{\pi-\theta}{\pi}, \quad \mbox{for any $\theta \in (0,\pi)$,}
\end{align*}
where $P_\theta$ is the convolution operator introduced in~\eqref{eq:convolutionPoisson}. Since the Poisson kernel $\mathscr{P}_\theta$ is positive, this estimate immediately follows from the fact that $\int_{\R} \mathscr{P}_\theta(\tau) \mathrm{d} \tau = \frac{\pi-\theta}{\pi}$.
\end{remark}

\begin{remark}[Simplified notation] To shorten the notation, from this point on we shall denote the composite map $\Psi_{\sigma_j}^\# \circ \gamma_j \circ (\Psi_{\sigma_i}^\#)^{-1}$ simply by $\gamma_j$. Note that, if we think of the pullback map $\Psi_{\sigma_j}^\#$ as a change of coordinates in the infinite-dimensional space $\mL^2(\partial \Omega_j)$, this slight abuse of notation is rather standard.
\end{remark}

\section{The adjoint Dirichlet-to-Dirichlet map in the Hardy space}\label{sec:duality}

In the previous section, we have shown (see Lemma~\ref{lem:DtDonStrip}) that, in order to prove Theorem~\ref{thm:DtDests}, it suffices to study the convolution operator
\begin{align*}
    (P_\theta g)(\tau) = \frac{1}{2\pi} \int_{\R} \frac{\sin(\theta)}{\cosh(\xi)-\cos(\theta)} g(\tau-\xi) \mathrm{d} \xi
\end{align*}
as an operator acting between the weighted spaces $\mathcal{L}^2(\R_\sigma)$ and $\mathcal{L}^2(\R_{\sigma\pm \theta})$. Therefore, we would like to apply the Fourier transform to turn $P_\theta$ into a multiplication operator. However, not all functions in $\mathcal{L}^2(\R_\sigma)$ have a well-defined Fourier transform as they may grow exponentially fast. To overcome this problem, let us introduce the following dual weighted spaces.
\begin{definition}[Dual weighted spaces] \label{def:dualweigthedspace} Let $1\leq p < \infty$ and $\sigma \in (0,2\pi)$, then we define $\mathcal{L}^{p^\ast}(\R_\sigma)$ as the vector space of measurable functions $f: \R \rightarrow \C$ such that
\begin{align*}
    \norm{f}_{\mathcal{L}^{p^\ast}(\R_\sigma)}^{\frac{p}{p-1}} \coloneqq \int_{\R} |f(\tau)|^{\frac{p}{p-1}}\left(\cosh(\tau) - \cos(\sigma)\right)  \mathrm{d} \tau < \infty.
\end{align*}
As before, we denote by $\mathcal{L}^{p^\ast}_\sigma$ the direct sum $\mathcal{L}^{p^\ast}(\R_\sigma) \oplus \mathcal{L}^{p^\ast}(\R_{\sigma+\pi}) = \mathcal{L}^{p^\ast}(\R_\sigma) \times \mathcal{L}^{p^\ast}(\R_{\sigma+\pi})$ endowed with the canonical norm 
\begin{align*}
    \norm{(f,g)}_{\mathcal{L}^{p^\ast}_\sigma}^{\frac{p}{p-1}} \coloneqq  \norm{f}_{\mathcal{L}^{p^\ast}(\R_\sigma)}^{\frac{p}{p-1}} + \norm{g}_{\mathcal{L}^{p^\ast}(\R_{\sigma+\pi})}^{\frac{p}{p-1}}.
\end{align*}
\end{definition}
Note that for $p=2$ the map 
\begin{align}
    W_\sigma : \mathcal{L}^2(\R_\sigma) \rightarrow \mL^2(\R) \quad\quad  f \mapsto (W_\sigma f) (\tau) = \frac{f(\tau)}{\sqrt{\cosh(\tau) - \cos(\sigma)}} \label{eq:weight}
\end{align}
is an isometric bijection. Moreover, the same map $W_\sigma$ also defines an isometric bijection from $\mL^2(\R)$ to $\mathcal{L}^{2^\ast}(\R_\sigma)$. More precisely, we have the following properties.
\begin{lemma}[Duality relations] \label{lem:dualityrelations} For any $\sigma \in (0,2\pi)$ the map $W_\sigma^2 : \mathcal{L}^2(\R_\sigma) \rightarrow \mathcal{L}^{2^\ast}(\R_\sigma)$ is an isometric bijection and we have
\begin{align*}
    \inner{W_\sigma^2 f, g}_{\mL^2(\R)} = \inner{f, g}_{\mathcal{L}^2(\R_\sigma)} = \inner{W_\sigma^2 f, W_\sigma^2 g}_{\mathcal{L}^{2^\ast}(\R_\sigma)} 
\end{align*}
for any $f,g \in \mathcal{L}^2(\R_\sigma)$. 
\end{lemma}
\begin{proof} The proof is immediate from the definition of the weigthed spaces and of $W_\sigma$.
\end{proof}

As a consequence of the above lemma, we have
\begin{lemma}[Adjoint Dirichlet-to-Dirichlet map] \label{lem:adjointDtD} Let $\mathbf{W}_\sigma : \mathcal{L}^2_{\sigma} \rightarrow \mL^2(\R)\times \mL^2(\R)$ be defined as
\begin{align*}
    \mathbf{W}_\sigma(f,g) \coloneqq (W_\sigma f, W_{\sigma+\pi}g).
\end{align*}
Then the adjoint of $\gamma_1 : \mathcal{L}^2_{\sigma_2} \rightarrow \mathcal{L}^2_{\sigma_1}$, i.e., the unique bounded operator $\gamma_1^\ast : \mathcal{L}^2_{\sigma_1} \rightarrow \mathcal{L}^2_{\sigma_2}$ such that
\begin{align*}
    \inner{\mathbf{f},\gamma_1 \mathbf{g}}_{\mathcal{L}^2_{\sigma_1}} = \inner{\gamma_1^\ast \mathbf{f}, \mathbf{g}}_{\mathcal{L}^2_{\sigma_2}} \quad \mbox{for any $\mathbf{f}\in \mathcal{L}^2_{\sigma_1}$ and any $\mathbf{g}\in \mathcal{L}^2_{\sigma_2}$,}
\end{align*}
satisfies $\mathbf{W}_{\sigma_2}^{2} \gamma_1^\ast  \mathbf{W}_{\sigma_1}^{-2} = T_1$, where $T_1 : \mathcal{L}^{2^\ast}_{\sigma_1} \rightarrow \mathcal{L}^{2^\ast}_{\sigma_2}$ is given by
\begin{align}
    T_1(f,g) = (P_{\pi-\theta} g, P_{\theta} g). \label{eq:T1def}
\end{align}
where $\theta = \sigma_2-\sigma_1$. Similarly, $\mathbf{W}_{\sigma_1}^2 \gamma_2^\ast \mathbf{W}_{\sigma_2}^{-2}= T_2$, where $T_2: \mathcal{L}^{2^\ast}_{\sigma_2} \rightarrow \mathcal{L}^{2^\ast}_{\sigma_1}$ is given by
\begin{align}
    T_2(f,g) = \left(P_{\theta} f, P_{\pi-\theta} f\right). \label{eq:T2def}
\end{align}
In particular, we have $\mathbf{W}_{\sigma_1}^2 (\gamma_1 \gamma_2)^\ast \mathbf{W}_{\sigma_1}^{-2} = T_2 T_1 = T_{21} : \mathcal{L}^{2^\ast}_{\sigma_1} \rightarrow \mathcal{L}^{2^\ast}_{\sigma_1}$, where 
\begin{align*}
    T_{21}(f,g) = \left(P_\theta P_{\pi-\theta} g, P_{\pi-\theta}^2 g\right).
\end{align*}
\end{lemma}

The advantage of working in the dual space is that we can identify its Fourier transform\footnote{See~\eqref{eq:conventionFT} for the convention used for the Fourier transform.} with the following Hardy space.
\begin{definition}[Hardy space] We denote by $\mathbb{H}^2(\mathcal{S}^\ast)$ the set of holomorphic functions on the (dual) strip
\begin{align*}
    \mathcal{S}^\ast \coloneqq \R \times (-1/2,1/2) \subset \C
\end{align*}
such that
\begin{align}
    \norm{h}_{\mathbb{H}^2(\mathcal{S}^\ast)}^2 = \sup_{-1/2<y<1/2} \int_\R |h_y(x)|^2 \mathrm{d} x < \infty,
\end{align}
where $h_y:\R \rightarrow \C$ denotes the restriction $h_y(x) = h(x+\ii y)$.
\end{definition}

We then have the following characterization:
\begin{lemma}[Fourier transform of dual space]\label{lem:Hardy} Let $\sigma \in (0,2\pi)$, then the Fourier transform
\begin{align}
    f\mapsto (\mathcal{F} f)(z) = \widehat{f}(z) = \int_\R \ee^{-\ii z \tau} f(\tau) \mathrm{d} \tau, \quad z\in \mathcal{S}^\ast \label{eq:conventionFT}
\end{align}
defines a continuous bijection from $\mathcal{L}^{2^\ast}(\R_\sigma)$ to $\mathbb{H}^2(\mathcal{S}^\ast)$. Moreover, this map is an isometric bijection if we endow $\mathbb{H}^2(\mathcal{S}^\ast)$ with the equivalent norm
\begin{align}
    \norm{h}_{\mathbb{H}^2_\sigma}^2 \coloneqq \frac{1}{2\pi} \int_\R \frac{|h_{1/2}(x)|^2 + |h_{-1/2}(x)|^2}{2} - \cos(\sigma) |h_0(x)|^2 \mathrm{d} x, \label{eq:equivalentnorm}
\end{align}
where $h_{\pm 1/2}$ denotes the boundary values\footnote{which are well-defined by the proof of the lemma.} of $h$ along the upper and lower boundary of the strip $\mathcal{S}^\ast$, and $h_0(x) = h(x)$.
\end{lemma}

\begin{proof}
Let $f \in \mathcal{L}^{2^\ast}(\R)$, then for any $z\in \mathcal{S}^\ast$ we have
\begin{align*}
   \left| \int_\R \ee^{-\ii z \tau} f(\tau) \mathrm{d} \tau\right|^2 \leq  \left(\int_\R \left(\cosh(\tau) - \cos(\sigma)\right) f(\tau)|^2 \mathrm{d} \tau \right)\left(\int_\R \frac{\ee^{2y \tau}}{\cosh(\tau) - \cos(\sigma) } \mathrm{d} \tau \right)\lesssim_{\sigma,y} \norm{f}_{\mathcal{L}^{2^\ast}_\sigma}^2 .
\end{align*}
Thus the Fourier transform $\widehat{f}$ is pointwise well-defined in $\mathcal{S}^\ast$. Moreover, a straightforward application of Fubini's and Morera's theorem implies that $\widehat{f}$ is holomorphic on $\mathcal{S}^\ast$. On the other hand, from Plancherel's theorem, for any $|y| < 1/2$ we have
\begin{align*}
    \norm{\widehat{f}_y}_{\mL^2(\R)}^2 = 2\pi \norm{f \ee^{y \tau}}_{\mL^2(\R)}^2 \leq \sup_{\tau\in \R} \left\{\frac{\ee^{2y\tau}}{\cosh(\tau) - \cos(\sigma)}\right\} \norm{f}_{\mathcal{L}^{2^\ast}_\sigma}^2 = \frac{\norm{f}_{\mathcal{L}^{2^\ast}_\sigma}^2}{1-\cos(\sigma)}.
\end{align*}
As the above bound is independent of $|y| < 1/2$, we conclude that $\widehat{f} \in\mathbb{H}^2(\mathcal{S}^\ast)$.

For the opposite implication, a detailed proof can be found in \cite[Section 2]{Cor24}. For the sake of completeness, we only sketch the main steps. Let $h \in \mathbb{H}^2(\mathcal{S}^\ast)$, then the first step consists in using Cauchy's integral theorem to show that for any $\phi \in C_c^\infty(\R)$, 
\begin{align*}
    \inner{h_{y_0}, \widehat{\phi}_{y_1}} = \inner{h_{y_0+y}, \widehat{\phi}_{y_1-y}} \quad \mbox{for any $y_0,y,y_1 \in \R$ such that $\ii y_0 \in \mathcal{S}^\ast$ and $\ii y_0+y \in \mathcal{S}^\ast$.}
\end{align*}
This is justified by using a well-known pointwise bound on $h$ (see, e.g. \cite[Lemma 2.4]{Cor24} or \cite[Lemma 11.3]{Mas09}) and the fact that $\widehat{\phi}$ is entire and decays fast (locally uniformly in $y$) as $|x| \ra \infty$. From the above formula and Plancherel's theorem, we conclude that 
\begin{align}
    \widecheck{h}_y(\tau) = \widecheck{h}_0(\tau) \ee^{y \tau} \in \mL^2(\R) \quad \mbox{for almost every $\tau \in \R$ and any $|y| < 1/2$.}  \label{eq:exponentialfourier}
\end{align}
We can now use the uniform $\mL^2(\R)$ bound on $\{\widehat{h}_y\}_{|y|<1/2}$ to infer the existence of a weak limit as $y \ra \pm 1/2$.  By formula~\eqref{eq:exponentialfourier}, the weak limit is unique and satisfies the formula $\widehat{h}_{\pm 1/2}(\tau) = \widehat{h}_0(\tau) \ee^{\pm \frac12 \tau} \in \mL^2(\R)$. In particular, $\widehat{h}$ is square integrable against $\cosh(\tau) - \cos(\sigma)$ and belongs to $\mathcal{L}^{2^\ast}_\sigma$ for any $\sigma \in (0,\pi)$.

The proof of formula~\eqref{eq:equivalentnorm} is straightforward from Plancherel's theorem and~\eqref{eq:exponentialfourier}.
\end{proof}

\begin{remark}[Notation for $\sigma$-dependent Hardy space] Throughout the rest of the paper we shall often denote the norm in~\eqref{eq:equivalentnorm} simply by $\norm{h}_{\sigma}$. Moreover, we denote by $\mathbb{H}^2_\sigma(\mathcal{S}^\ast)$, or simply $\mathbb{H}^2_\sigma$, the Hardy space $\mathbb{H}^2(\mathcal{S}^\ast)$ endowed with this norm.
\end{remark}
We can now combine the previous results with the explicit formula\footnote{See, e.g., \cite[Lemma B.2]{Cor24} for a derivation of this formula.} for the Fourier transform of the Poisson kernel,
\begin{align}
    \widehat{\mathscr{P}}_\theta(z) = \frac{\sinh\left((\pi-\theta)z\right)}{\sinh(\pi z)}, \label{eq:PoissonFourier}
\end{align}
to obtain the following result, which is the main result of this section.

\begin{theorem}[Dirichlet-to-Dirichlet map in the Hardy space] \label{thm:DtDHardy} Let $\gamma_j$ be the DtD map introduced in Definition~\eqref{def:trace}. Then there exists $\sigma_1 < \sigma_2 \in (0,\pi)$ such that the operators $\gamma_1^\ast$, $\gamma_2^\ast$ and $(\gamma_1 \gamma_2)^\ast$ are, respectively, unitarily equivalent to the operators
\begin{align}
    &T_1 : \mathbb{H}^2_{\sigma_1} \times \mathbb{H}^2_{\sigma_1+\pi} \rightarrow \mathbb{H}^2_{\sigma_2} \times \mathbb{H}^2_{\sigma_2+\pi}, \quad T_1(f,g) = \left(\widehat{\mathscr{P}}_{\pi-\theta} g\,,\, \widehat{\mathscr{P}}_{\theta}g\right)\label{eq:T1Hardydef} \\
     &T_2 : \mathbb{H}^2_{\sigma_2} \times \mathbb{H}^2_{\sigma_2+\pi} \rightarrow \mathbb{H}^2_{\sigma_1} \times \mathbb{H}^2_{\sigma_1+\pi}, \quad T_2(f,g) = \left(\widehat{\mathscr{P}}_{\theta} f\,,\, \widehat{\mathscr{P}}_{\pi-\theta}f \right) \label{eq:T2Hardydef} \\
      &T_{21} : \mathbb{H}^2_{\sigma_1} \times \mathbb{H}^2_{\sigma_1+\pi} \rightarrow \mathbb{H}^2_{\sigma_1} \times \mathbb{H}^2_{\sigma_1+\pi}, \quad T_{21}(f,g) = \left(\widehat{\mathscr{P}}_\theta \widehat{\mathscr{P}}_{\pi-\theta} g\,,\,  \widehat{\mathscr{P}}_{\pi-\theta}^2 g\right),\label{eq:T21Hardydef}
\end{align}
where $\theta = \sigma_2 - \sigma_1$ and $\widehat{\mathscr{P}}_\theta$ is the Fourier transform of the Poisson kernel given in~\eqref{eq:PoissonFourier}.
\end{theorem}

\begin{proof} The proof is immediate from the adjoint operator formulas in Lemma~\ref{lem:adjointDtD}, the identification of the Fourier transform of $\mathcal{L}^{2^\ast}_\sigma$ with $\mathbb{H}^2(\mathcal{S}^\ast)$ in Lemma~\ref{lem:Hardy}, and the convolution property of the Fourier transform.
\end{proof}

\section{Spectral analysis of the Dirichlet-to-Dirichlet map}
\label{sec:proof}

According to Theorem~\ref{thm:DtDHardy}, all the statements in Theorem~\ref{thm:DtDests}, which concern the operators $\gamma_1$, $\gamma_2$ and $\gamma_1\gamma_2$, can instead be proven for the operators $T_1$, $T_2$, and $T_{21}$ introduced in~\eqref{eq:T1Hardydef}--\eqref{eq:T21Hardydef}. The goal of this section is therefore to prove these corresponding statements one by one. To this end, we shall use a few technical lemmas, whose proofs will be postponed to the appendix.

The first lemma contains some useful properties of the Fourier transform of the Poisson kernel.

\begin{lemma}[Properties of Poisson kernel] \label{lem:Poissonproperties} Let $\theta \in (0,\pi)$ and $\widehat{\mathscr{P}}_\theta$ be the Fourier transform of the Poisson kernel given by~\eqref{eq:PoissonFourier}. Then the following holds:
\begin{enumerate}[label=(\roman*)]
\item \label{it:Poissonmaxima} (Global maxima) The maximum of $\widehat{\mathscr{P}}_\theta$ over $\overline{\mathcal{S}^\ast}$ is achieved only at $\pm \ii/2$ and given by
\begin{align}
    \max_{z \in \overline{\mathcal{S}^\ast}} |\widehat{\mathscr{P}}_{\theta}(z)| = |\widehat{\mathscr{P}}_\theta(\pm \ii/2)| = \cos\left(\theta/2\right).
\end{align}
\item \label{it:combinedsupremumnorm} (Bound on difference along $\R$) For any $\theta \in (0,\pi/2]$ we have
\begin{align*}
    0 \leq |\widehat{\mathscr{P}}_\theta(x)|^2 - |\widehat{\mathscr{P}}_{\pi-\theta}(x)|^2 \leq 1-\frac{2\theta}{\pi} \leq \cos(\theta) \quad \mbox{for any $x\in \R$.} 
\end{align*}
\item (Reproducing kernel) \label{it:reproducingkernel} Let $\sigma \in (0,\pi)$ and $\inner{\cdot,\cdot}_{\sigma}$ denote the inner-product associated to the norm~\eqref{eq:equivalentnorm}. Then for any $h\in \mathbb{H}^2(\mathcal{S}^\ast)$ and $z_0 \in \mathcal{S}^\ast$, we have
\begin{align*}
    h(z_0) = \frac{2\pi}{\sin(\sigma)} \inner{\mathscr{T}_{z_0}\widehat{\mathscr{P}}_{\sigma},h}_\sigma,
\end{align*}
where $\mathscr{T}_{z_0} \widehat{\mathscr{P}}_\sigma$ denotes the translated map $z \mapsto \mathscr{T}_{z_0} \widehat{\mathscr{P}}_\sigma(z) \coloneqq \widehat{\mathscr{P}}_\sigma(z-\overline{z_0})$.
\end{enumerate}
\end{lemma}

The second lemma is a version of the residue theorem suited for our purposes. 
\begin{lemma}[Residue theorem in the Hardy space] \label{lem:residuethm} Let $F: \overline{\mathcal{S}^\ast}\setminus \{0\} \rightarrow \C$ be a function such that
\begin{enumerate}[label=(\roman*)]
\item $F$ is holomorphic in $\mathcal{S}^\ast\setminus \{0\}$,
\item $F$ is continuous up to the boundary and bounded away from $z=0$, and 
\item $F$ has at most a simple pole at $0$ with purely imaginary residue, i.e.
\begin{align*}
    \mathrm{Re} \{ \mathrm{Res} (F,0) \} = \mathrm{Re} \{\lim_{z \ra 0} F(z) z\} = 0.
\end{align*}
\end{enumerate}
Then $\mathrm{Re} \{F(x)\}$ is bounded for $x\in \R$ and we have
\begin{align}
 \mathrm{Re}\left\{\int_\R F(x+\ii/2)  h_{1/2}(x) \overline{h_{-1/2}(x)}  \mathrm{d} x\right\} = \int_\R \mathrm{Re} \{ F(x)\} |h_0(x)|^2 \mathrm{dx} - \ii \pi \mathrm{Res}(F,0) |h(0)|^2, \label{eq:residuethm}
\end{align}
for any $h \in \mathbb{H}^2(\mathcal{S}^\ast)$. 
\end{lemma}

The third lemma is a localization procedure in $\mathbb{H}^2(\mathcal{S}^\ast)$ that will be useful to prove the lower bounds in Theorem~\ref{thm:DtDests}.
\begin{lemma}[Localization on Hardy space] \label{lem:localizationlemma} Let $f \in \mL^2(\R)$, then the function 
\begin{align}
    (m f)(z) \coloneqq \int_{\R} 2\widehat{\mathscr{P}}_{\pi/2}(2z-2u) f(u) \mathrm{d} u \label{eq:moperator}
\end{align}
belongs to $\mathbb{H}^2(\mathcal{S}^\ast)$ and satisfies
\begin{align}
    \frac{(mf)_{1/2}(x) + (mf)_{-1/2}(x)}{2} = f(x)\quad \mbox{for almost every $x\in \R$.} \label{eq:firsttransformidentity}
\end{align}
Moreover, if we define the dilation operator $D_\delta : \mL^2(\R) \rightarrow \mL^2(\R)$ as $(D_\delta f)(x) \coloneqq f(x/\delta)\delta^{-1/2}$ for any $\delta >0$, then we have
\begin{align}
    \lim_{\delta \downarrow 0} (m D_\delta f)_0 = 0 \quad \mbox{and}\quad \lim_{\delta \downarrow 0} D_{1/\delta} \left(\frac{(m D_\delta f)_{1/2}- (m D_\delta f)_{-1/2}}{2}\right) = -\ii H(f) \label{eq:secondtransformidentity}
\end{align}
where the limit is taken in the $\mL^2(\R)$ sense, and $H(f)$ is the Hilbert transform of $f$,
\begin{align*}
    (Hf)(x) = \lim_{\epsilon \downarrow 0^+} \int_{|y| \geq \epsilon} \frac{f(x-y)}{\pi y} \mathrm{d} y.
\end{align*}
\end{lemma}

Equipped with these lemmas, we can now proceed to the proof of Theorem~\ref{thm:DtDests}.

\subsection{Proof of Theorem~\ref{thm:DtDests}~\ref{it:operatornorm}}
We start with the proof of the following theorem, which is an improved version of Theorem~\ref{thm:DtDests}~\ref{it:operatornorm}.

\begin{theorem}[General bound on $\gamma_2$] \label{thm:generalbound} Let $\gamma_2 : \mathrm{L}^2(\partial \Omega_1) \rightarrow \mL^2(\partial \Omega_2)$ be the DtD map introduced in Definition~\ref{def:trace} and let $\theta,\sigma \in(0,\pi)$ be the angles described in Section~\ref{sec:bipolartransform}. Then we have
\begin{align}
    1\leq \norm{\gamma_2}^2 \leq 1+ \left( \frac{\sin(\theta)}{\tan(\sigma+\theta)}\right)_+ \frac{\pi-\sigma-\theta}{\pi}. \label{eq:detailedoperatornorm}
\end{align}
Moreover, there exists $\theta,\sigma \in (0,\pi)$ such that $\norm{\gamma_2} > 1$. 
\end{theorem}

\begin{remark*}[Sharp norm] Recalling the convention adopted in Section~\ref{sec:bipolartransform}, we note that $\pi-\beta_2 = \sigma+\theta$, where $\beta_2$ is the angle in Figure~\ref{fig:main_def}. In particular, estimate~\eqref{eq:detailedoperatornorm} yields the sharp bound $\norm{\gamma_2} =1$ for $\beta_2 \leq \pi/2$. 
\end{remark*}
\begin{proof} By Theorem~\ref{thm:DtDHardy} and the well-known identity $\norm{\gamma_2} = \snorm{\gamma_2^\ast}$, it suffices to prove~\eqref{eq:detailedoperatornorm} for the operator $T_2 : \mathbb{H}^2_{\sigma+\theta} \rightarrow \mathbb{H}^2_{\sigma}\times \mathbb{H}^2_{\sigma+\pi}$ defined as $T_2(h) = (\widehat{\mathscr{P}}_\theta h, \widehat{\mathscr{P}}_{\pi-\theta} h)$ for any $\theta,\sigma \in (0,\pi)$ with $\theta + \sigma \leq \pi$. So let us start with the upper bound. To shorten the notation throughout the proof, let us define
\begin{align*}
    F_\theta(x) \coloneqq |\widehat{\mathscr{P}}_\theta(x+\ii/2)|^2 + |\widehat{\mathscr{P}}_{\pi-\theta}(x+\ii/2)|^2 \quad\mbox{and}\quad G_\theta(x) \coloneqq |\widehat{\mathscr{P}}_\theta(x)|^2 - |\widehat{\mathscr{P}}_{\pi-\theta}(x)|^2. 
\end{align*}
Then, from the symmetry $|\widehat{\mathscr{P}}_\theta(x+\ii/2)| = |\widehat{\mathscr{P}}_\theta(x-\ii/2)|$, we note that
\begin{align}
    \norm{h}_{\mathbb{H}^2_{\sigma+\theta}}^2 - \norm{T_2 h}_{\mathbb{H}^2_{\sigma} \times \mathbb{H}^2_{\sigma+\pi}}^2 &= \frac{1}{2\pi} \int_{\R} \left(1-F_\theta(x)\right) \frac{|h_{1/2}(x)|^2 + |h_{-1/2}(x)|^2}{2}\mathrm{d}x \nonumber \\
    &+ \frac{1}{2\pi} \int_\R \left(\cos(\sigma)G_\theta(x)-\cos(\sigma+\theta)\right)|h_0(x)|^2 \mathrm{d}x \label{eq:normexpression}
\end{align}
Next, we apply the residue theorem in Lemma~\ref{lem:residuethm} to the function 
\begin{align*}
    F(z) \coloneqq 1- \widehat{\mathscr{P}}_\theta(z)\widehat{\mathscr{P}}_\theta(z-\ii) - \widehat{\mathscr{P}}_{\pi-\theta}(z) \widehat{\mathscr{P}}_{\pi-\theta}(z-\ii) = 1-F_\theta(z-\ii/2).
\end{align*}
Precisely, we note that
\begin{align*}
    \mathrm{Re}\{ F(x)\} = 1-\cos(\theta) G_\theta(x) \quad \mbox{and}\quad \mathrm{Res}(F,0) = -\ii \frac{\sin(\theta)}{\pi},
\end{align*}
and therefore, by the pointwise identity
\begin{align}
    |h_{1/2}|^2 + |h_{-1/2}|^2 = (1-\alpha) \frac{|h_{1/2}|^2 + |h_{-1/2}|^2}{2} + \alpha \frac{|h_{1/2} - h_{-1/2}|^2}{2} + \alpha \mathrm{Re} \,\{ h_{1/2} \overline{h_{-1/2}}\},\label{eq:pointwiseid}
\end{align}
for $\alpha \in \R$, and the residue theorem in~\eqref{eq:residuethm}, the right-hand side of~\eqref{eq:normexpression} is equal to
\begin{align*}
    \frac{1-\alpha}{2\pi}\int_\R \left(1-F_\theta(x)\right) \frac{|h_{1/2}(x)|^2 + |h_{-1/2}(x)|^2}{2} + \frac{\alpha}{2\pi} \int_\R \left(1-F_\theta(x)\right) \frac{|h_{1/2}(x)-h_{-1/2}(x)|^2}{2} \mathrm{d} x\\
    + \frac{1}{2\pi} \int_\R \left(\alpha - \cos(\sigma+\theta) + \left(\cos(\sigma) - \alpha \cos(\theta)\right) G_\theta(x)\right) |h_0(x)|^2 \mathrm{d} x - \frac{\alpha \sin(\theta)}{2\pi} |h(0)|^2.
\end{align*}
As $F_\theta(x) \leq 1$ by Lemma~\ref{lem:Poissonproperties}~\ref{it:Poissonmaxima}, the first two terms are non-negative as long as $\alpha \in [0,1]$. The idea now is to find the smallest possible $\alpha \in [0,1]$ such that the integrand in the third integral is also non-negative.

We claim that this $\alpha$ is given by
\begin{align*}
    \alpha = \cos(\sigma+\theta)_+.
\end{align*}
To prove this claim in the case $\theta+\sigma \leq \pi/2$, we first note that $\alpha$ must be at least $\cos(\sigma+\theta)$ because
\begin{align*}
    \lim_{|x| \ra \infty} \alpha - \cos(\sigma+\theta) + \left(\cos(\sigma) - \alpha \cos(\theta)\right) G_\theta(x) = \alpha-\cos(\sigma+\theta).
\end{align*}
Moreover, to see that $\cos(\theta+\sigma)$ suffices, we note that $\theta \leq \pi/2$ and therefore $G_\theta(x) \geq 0$; hence,
\begin{align*}
    \left(\cos(\sigma) - \cos(\sigma+\theta) \cos(\theta)\right) G_\theta(x) = \sin(\sigma+\theta) \sin(\theta) G_\theta(x) \geq 0,
\end{align*}
and the claim follows. To prove the claim in the case $\sigma+\theta \geq \pi/2$, we need to show that
\begin{align*}
    -\cos(\sigma+\theta) + \cos(\sigma) G_\theta(x) \geq 0\quad \mbox{for any $x\in \R$ and $\sigma+\theta \geq \pi/2$.}
\end{align*}
To this end, notice that, due to the restriction $\theta+\sigma \leq \pi$, we can not have simultaneously $\sigma \geq \pi/2$ and $\theta \geq \pi/2$. As $G_\theta(x)$ has the same sign as $\cos(\theta)$, this implies that $G_\theta(x)$ and $\cos(\sigma)$ can not be simultaneously negative. Since, in addition, $|G_\theta(x)| \leq |\cos(\theta)|$ by Lemma~\ref{lem:Poissonproperties}~\ref{it:combinedsupremumnorm}, we have
\begin{align*}
    -\cos(\sigma+\theta) + \cos(\sigma) G_\theta(x) \geq -\cos(\sigma+\theta) - \left(-\cos(\sigma) \cos(\theta)\right)_+. 
\end{align*}
Thus, on one hand, if $\cos(\theta)\cos(\sigma) \geq 0$, the claim holds because $-\cos(\sigma+\theta) \geq 0$ for $\theta +\sigma \in (\pi/2,\pi)$. On the other hand, if $\cos(\theta) \cos(\sigma) \leq 0$, the claim also holds because we have $-\cos(\sigma+\theta) +\cos(\sigma) \cos(\theta) = \sin(\theta) \sin(\sigma) \geq 0$ as both $\theta,\sigma \in (0,\pi)$.

We have thus shown that
\begin{align}
    \norm{h}_{\sigma+\theta}^2 - \norm{T_2 h}_{\sigma,\sigma+\pi}^2 \geq -\cos(\theta+\sigma)_+  \frac{\sin(\theta)}{2\pi} |h(0)|^2. \label{eq:somelowerbound}
\end{align}
Hence, we can now use Lemma~\ref{lem:Poissonproperties}~\ref{it:reproducingkernel} with $z_0 = 0$ and Cauchy-Schwarz to obtain
\begin{align}
    -\frac{\cos(\sigma+\theta) \sin(\theta)}{2\pi} |h(0)|^2 &= -\cos(\sigma+\theta)_+ \sin(\theta) \frac{2\pi}{\sin(\theta+\sigma)^2}|\inner{\mathscr{P}_{\sigma+\theta},h}|^2 \nonumber \\
    &\geq -\cos(\sigma+\theta)_+ \sin(\theta) \frac{2\pi}{\sin(\theta+\sigma)^2}\inner{\widehat{\mathscr{P}}_{\sigma+\theta}, \widehat{\mathscr{P}}_{\sigma+\theta}}_{\sigma+\theta} \norm{h}_{\sigma+\theta}^2 \nonumber \\
    &= \left(\frac{\sin(\theta)}{\tan (\sigma+\theta)}\right)_+ \frac{\pi-\sigma-\theta}{\pi} \norm{h}_{\sigma+\theta}^2, \label{eq:CS}
\end{align}
where the last equality follows from Lemma~\ref{lem:Poissonproperties}~\ref{it:reproducingkernel} since $\widehat{\mathscr{P}}_{\sigma+\theta}(0) = (\pi-\sigma-\theta)/\pi$. This bound together with~\eqref{eq:somelowerbound} proves the upper bound in~\eqref{eq:detailedoperatornorm}.

To prove the uniform lower bound $1\leq \norm{\gamma_2}$, we use the localization procedure in Lemma~\ref{lem:localizationlemma}. So let $f\in \mL^2(\R)$ and define the trial state
\begin{align*}
    h^\delta(z) \coloneqq \left(m D_\delta f\right)(z), 
\end{align*}
where $m$ is defined according to~\eqref{eq:moperator}. Then, on one hand, from the identity
\begin{align*}
    \norm{h^\delta}_{\sigma+\theta}^2 &= \frac{1}{2\pi} \int_{\R} \frac{|h_{1/2}^\delta(x)|^2 + |h_{-1/2}^\delta(x)|^2}{2} - \cos(\sigma+\theta) |h^\delta_0(x)|^2 \mathrm{d} x \\
    &= \frac{1}{2\pi} \norm{\frac{h_{1/2}^\delta+h_{1/2}^\delta}{2}}_{\mL^2(\R)}^2 + \frac{1}{2\pi} \norm{\frac{h_{1/2}^\delta-h_{-1/2}^\delta}{2}}_{\mL^2(\R)}^2 - \frac{\cos(\sigma+\theta)}{2\pi} \norm{h^\delta_0}_{\mL^2(\R)}^2
\end{align*}
and properties~\eqref{eq:firsttransformidentity} and~\eqref{eq:secondtransformidentity} of the $m$ operator, we find that
\begin{align}
    \lim_{\delta \downarrow 0} \norm{h^\delta}_{\mathbb{H}^2_{\sigma+\theta}}^2 = \frac{1}{2\pi} \norm{f}_{\mL^2(\R)}^2 + \frac{1}{2\pi}\norm{H f}_{\mL^2(\R)}^2 = \frac{1}{\pi} \norm{ f}_{\mL^2(\R)}^2. \label{eq:hdeltanorm}
\end{align}
On the other hand, by re-scaling, we have
\begin{align*}
    \norm{T_2 h^\delta}&= \frac{1}{2\pi} \int_\R F_\theta(x) \frac{|h^\delta_{1/2}(x)|^2 + |h^\delta_{-1/2}(x)|^2}{2} -\cos(\sigma) G_\theta(x) |h^\delta(x)|^2 \mathrm{d} x\\
    &= \frac{1}{2\pi} \int_\R F_\theta(x/\delta) \left(|f(x)|^2 + \left|D_{1/\delta}\left(\frac{h^{\delta}_{1/2} - h_{-1/2}^\delta}{2}\right)(x)\right|^2\right)\mathrm{d}x - \frac{\cos(\sigma)}{2\pi} \int_\R  G_\theta(x)|h^\delta(x)|^2 \mathrm{d} x.
\end{align*}
Since $F_\theta(0) = 1$ and both $F_\theta$ and $G_\theta$ are uniformly bounded and continuous in $\R$, we can use properties~\eqref{eq:firsttransformidentity} and~\eqref{eq:secondtransformidentity} and the dominated convergence theorem to conclude that
\begin{align*}
    \lim_{\delta \downarrow 0} \norm{T_2 h^\delta}_{\mathbb{H}^2_\sigma\times \mathbb{H}^2_{\sigma+\pi}}^2 = \frac{1}{\pi} \norm{f}^2_{\mL^2(\R)}, 
\end{align*}
which together with~\eqref{eq:hdeltanorm} shows that $\norm{\gamma_2} \geq 1$. 

To prove the strict lower bound $\norm{\gamma_2}>1$ for certain values of $\theta$ and $\sigma$, we shall argue by contradiction. So suppose that $\norm{\gamma_2}\leq 1$ for any $\theta,\sigma \in (0,\pi)$ with  $\theta+\sigma \leq \pi$. Then 
\begin{align*}
    E(h,\theta) &\coloneqq  \frac{1}{2\pi} \int_{\R} \left(1-F_\theta(x)\right) \frac{|h_{1/2}(x)|^2+|h_{-1/2}(x)|^2}{2} \mathrm{d} x + \frac{1}{2\pi} \int_\R \left(G_\theta(x) - \cos(\theta)\right) |h_0(x)|^2 \mathrm{d} x \\
    &= \lim_{\sigma \downarrow 0}  \left(\norm{h}_{\mathbb{H}^2_{\sigma+\theta}}^2 - \norm{T_2 h}_{\mathbb{H}^2_\sigma \times \mathbb{H}^2_{\sigma+\pi}}^2\right) \geq 0, \quad \mbox{for any $\theta \in (0,\pi)$ and $h \in \mathbb{H}^2(\mathcal{S}^\ast)$.}
\end{align*}
Next, we consider the limit $\theta \downarrow 0$. For this, we first note that
\begin{align*}
    |1-F_\theta(x)| \leq \theta|x| \quad \mbox{and}\quad |G_\theta(x) - \cos(\theta)| \leq \theta |x|\quad  \mbox{uniformly in $\theta \in (0,\pi)$ and $x \in \R$.}
\end{align*}
Moreover, from straightforward calculation, we find
\begin{align*}
    \lim_{\theta \downarrow 0} \frac{1-F_\theta(x)}{\theta} = 2x\tanh(\pi x) \quad \mbox{and}\quad \lim_{\theta \downarrow 0} \frac{G_\theta(x) - \cos(\theta)}{\theta} = -2 x\coth(\pi x).
\end{align*}
Therefore, for any $h \in \mathbb{H}^2(\mathcal{S}^\ast)$ satisfying
\begin{align*}
    \sup_{-1/2 \leq y \leq 1/2}\int_\R (1+|x|) |h_y(x)|^2 \mathrm{d} x < \infty,
\end{align*}
 we can apply the dominated convergence theorem to obtain
 \begin{align*}
     \partial_\theta E(h,0) &= \lim_{\delta \downarrow 0} \frac{E(h,\theta)}{\theta} \\
     &= \frac{1}{2\pi} \int_\R 2x \tanh(\pi x) \frac{|h_{1/2}(x)|^2 + |h_{-1/2}(x)|^2}{2} \mathrm{d} x - \frac{1}{2\pi} \int_\R \frac{2x}{\tanh(\pi x)} |h_0(x)|^2 \mathrm{d} x 
 \end{align*}
We now use the residue theorem to cancel the second integral. Precisely, we note that
\begin{align*}
    2(x-\ii/2) \tanh(\pi (x-\ii/2) = (2x-\ii) \coth(\pi x),
\end{align*}
and therefore, an application of the residue theorem in~\eqref{eq:residuethm} shows that
 \begin{align*}
     \partial_\theta E(h,0) = \frac{1}{2\pi} \int_\R 2x \tanh(\pi x) \frac{|h_{1/2}(x)-h_{-1/2}(x)|^2}{2} \mathrm{d} x - \frac{|h(0)|^2}{2\pi}
 \end{align*}
 Let us now consider trial states of the form $h^\delta(z) = h(\delta z)$ for $0<\delta <1$ and $h \in \mathbb{H}^2(\mathcal{S}^\ast)$ satisfying 
\begin{align}
    \sup_{-1/2<y<1/2} \int_\R (1+|x|)|\partial_y h_y(x)|^2 \mathrm{d} x < \infty. 
\end{align}
For instance, any $h = \widehat{\phi}$ with $\phi \in C_c^\infty(\R)$ is fine. Then, since
\begin{align*}
    |h^\delta_{1/2}(x) - h^\delta_{-1/2}(x)|^2 &= |h_{\delta/2}(\delta x) - h_{-\delta/2}(\delta x)|^2 \\
    &= \delta^2 \int_{-1}^1 \int_{-1}^1 \int_{-1}^1 \partial_y h_{-\delta/2 + t\delta}(\delta x) \overline{\partial_y h_{-\delta/2 + s\delta}(\delta x)} \mathrm{d}s \mathrm{d} t,
\end{align*}
 by re-scaling in $x$ and applying dominated convergence we find
\begin{align*}
    q(h) &\coloneqq \lim_{\delta \downarrow 0^+} \partial_\theta E(h^\delta,0) \\
    &= \lim_{\delta \downarrow 0^+} \frac{1}{2\pi} \int_{-1}^1\int_{-1}^1 \int_\R  x \tanh(\pi x) \delta^2 \partial_y h_{-\delta/2+t\delta}(\delta x) \overline{\partial_y h_{-\delta/2+s\delta}(x)} \mathrm{d} x\mathrm{d}s \mathrm{d} t - \frac{|h(0)|^2}{\pi} \\
    &= \frac{2}{\pi} \int_\R |x| |\partial_x h(x)|^2 \mathrm{d} x - \frac{|h(0)|^2}{2\pi}
\end{align*}
where we used in the last line that $|\partial_y h(z)| = |\partial_x h(z)|$ by the Cauchy-Riemann equations. As we assumed $E(h,\theta) \geq 0$ for any $h\in \mathbb{H}^2(\mathcal{S}^\ast)$ and $\theta \in (0,\pi)$, all the limits considered so far are limits of non-negative numbers. In particular, we must have $q(h) \geq 0$ for any $h = \widehat{\phi} \in C^\infty_c(\R)$. Hence, to conclude the proof, it suffices to show that $q(h)<0$ for some $h = \widehat{\phi}$ with $\phi \in C_c^\infty(\R)$. For this, let us define the function
\begin{align*}
    f(x) \coloneqq \begin{dcases} \log(K), \quad &\mbox{ for $|x| \leq 1$}, \\
    \log(K) -\log(|x|), \quad &\mbox{for $1 \leq |x| \leq K$}, \\
    0, \quad &\mbox{otherwise.} 
    \end{dcases}
\end{align*}
and note that
\begin{align*}
    q(f) = \frac{2}{\pi} \int_\R |x| |\partial_x f(x)|^2 \mathrm{d} x - \frac{\log (K)^2}{2\pi} = 4\int_1^K x^{-1} \mathrm{d}x - \frac{\log(K)^2}{2\pi} = \frac{4}{\pi} \log(K)  -\frac{\log(K)^2}{2\pi} < 0 , 
\end{align*}
for $K>0$ large enough. Clearly $f$ is not an admissible trial state as $\widehat{f} \not \in C^\infty_c(\R)$. However, since $f$ has compact support, we must have $\widehat{f} \in C^\infty(\R)$, and therefore $\mathcal{F}(f \ast \phi) = \widehat{f} \widehat{\phi} \in C_c^\infty(\R)$ for any $\phi$ such that $\widehat{\phi} \in C^\infty_c(\R)$. In particular, the sequence of approximating functions $h^\epsilon = f \ast \phi^\epsilon$, where $\phi^\epsilon$ are mollifiers with $\widehat{\phi}^\epsilon \in C^\infty_c(\R)$, are admissible and satisfy
\begin{align*}
    \lim_{\epsilon \downarrow 0} q(h^\epsilon) = q(h) < 0.
\end{align*}
Hence $q(h^\epsilon)<0$ for $\epsilon$ small, which yields the desired contradiction and completes the proof.
\end{proof}

\subsection{Proof of Theorem~\ref{thm:DtDests}~\ref{it:spectralradius}} We now turn to the proof of the following result, which combined with Lemma~\ref{lem:Poissonproperties}~\ref{it:Poissonmaxima} immediately implies Theorem~\ref{thm:DtDests}~\ref{it:spectralradius}. 

\begin{theorem}[Spectrum of $\gamma_1 \gamma_2$] \label{thm:spectrum} The spectrum of the operator $\gamma_1 \gamma_2$ is given by
\begin{align*}
    \sigma(\gamma_1 \gamma_2) = \mathrm{cl} \left\{ \overline{\lambda} \in \C : \lambda \in \mathscr{P}^2_{\pi-\theta}(\mathcal{S}^\ast)\right\},
\end{align*}
where $\mathrm{cl}\, U$ denotes the closure of $U\subset \C$. In particular, $\rho(\gamma_1 \gamma_2) = \frac{1-\cos(\theta)}{2}$. Moreover, any point $\overline{\lambda} = \mathscr{P}^2_{\pi-\theta}(z_0)$ for $z_0 \in \mathcal{S}^\ast$ is an eigenvalue of $\gamma_1 \gamma_2$ with eigenfunction
\begin{align}
    f(\phi) = \left(\frac{1-\cos(\phi+\sigma)}{\cos(\phi)-\cos(\sigma)}\right)^{\ii \overline{z_0}} \in \mL^2(\Gamma^{\rm int}_1). \label{eq:expliciteigenfunction}
\end{align}
\end{theorem}
\begin{proof}
    From the identity $\sigma(\gamma_1\gamma_2) = \{z \in \C: \overline{z} \in \sigma\left((\gamma_1 \gamma_2)^\ast\right)\}$ and Theorem~\ref{thm:DtDHardy}, to prove the first statement it suffices to show that the spectrum of the operator 
    \begin{align*}
    T_{21}(g,h) = \left(\widehat{\mathscr{P}}_\theta \widehat{\mathscr{P}}_{\pi-\theta} g\,,\,  \widehat{\mathscr{P}}_{\pi-\theta}^2 g\right)
    \end{align*}
    is given by the closure of the image of the strip $\mathcal{S}^\ast$ under the map $z\mapsto \widehat{\mathscr{P}}_{\pi-\theta}(z)^2$. To show this, the key observation is that the inverse of $\lambda- T_{21}$ can be explicitly computed pointwise, i.e., given $(f,g) \in \mathbb{H}^2(\mathcal{S}^\ast)$, then 
    \begin{align}
        (\lambda - T_{21})(\tilde{f},\tilde{g}) = \left(\lambda \tilde{f}- \widehat{\mathscr{P}}_\theta \widehat{\mathscr{P}}_{\pi-\theta} \tilde{g}, \left(\lambda - \widehat{\mathscr{P}}_{\pi-\theta}^2\right) \tilde{g}\right) =  (f,g) \label{eq:solution}
    \end{align}
    for functions $\tilde{f}, \tilde{g} : \mathcal{S}^\ast \rightarrow \C$ if and only if
    \begin{align}
        \tilde{f}(z) = \frac{1}{\lambda}\left( f(z) + \frac{\widehat{\mathscr{P}}_{\pi-\theta}(z) \widehat{\mathscr{P}}_\theta(z)}{\lambda - \widehat{\mathscr{P}}_{\pi-\theta}(z)^2} g(z)\right) \quad \mbox{and} \quad \tilde{g}(z) = \frac{1}{\lambda - \mathscr{P}_{\pi-\theta}(z)^2} g(z), \quad \mbox{for $z\in \mathcal{S}^\ast$.} \label{eq:inverseT21}
    \end{align}
    Hence, if $\lambda = \widehat{\mathscr{P}}_{\pi-\theta}(z_0)^2$ for some $z_0\in \mathcal{S}^\ast$, the function
    \begin{align*}
        z\mapsto (\lambda - \widehat{\mathscr{P}}_{\pi-\theta}(z)^2)^{-1} g(z)
    \end{align*}
    has a pole at $z_0$ for any  $g\in \mathbb{H}^2(\mathcal{S}^\ast)$ with $g(z_0) \neq 0$ (e.g. $g(z) = (z+i)^{-2}$) and there is no solution of~\eqref{eq:solution} in $\mathbb{H}^2(\mathcal{S}^\ast)^2$. In particular, $\lambda - T_{21}$ is not surjective and $\lambda \in \sigma(T_{21})$. As the spectrum is a closed set, this implies that
    \begin{align*}
        \mathrm{cl} \, \widehat{\mathscr{P}}_{\pi-\theta}^2(\mathcal{S}^\ast) \subset \sigma(T_{21}).
    \end{align*}
    
    To prove the opposite inclusion, we note that multiplication by a bounded holomorphic function is a bounded operator in $\mathbb{H}^2(\mathcal{S}^\ast)$. Indeed, since continuity is independent of the choice of equivalent norm, we can consider the original norm in $\mathbb{H}^2(\mathcal{S}^\ast)$ and use the trivial inequality
    \begin{align*}
        \norm{f h}_{\mathbb{H}^2(\mathcal{S}^\ast)}^2 = \sup_{-1/2<y<1/2} \norm{ f_y h_y}_{\mL^2(\R)}^2 \leq \norm{f}_{\mL^\infty(\mathcal{S}^\ast)}^2 \sup_{-1/2<y<1/2} \norm{ h_y}_{\mL^2(\R)}^2,
    \end{align*}
    to conclude that $h \mapsto f h$ is continuous in $\mathbb{H}^2(\mathcal{S}^\ast)$ for any bounded holmorphic function $f$. Since the function $1/(\lambda - \widehat{\mathscr{P}}_{\pi-\theta}^2)$ is uniformly bounded in $\mathcal{S}^\ast$ as long as $\lambda$ is outside the closure of the image of $\widehat{\mathscr{P}}_{\pi-\theta}^2$, the function $\tilde{g}$ defined in~\eqref{eq:solution} is holomorphic and satisfies
    \begin{align*}
        \norm{\tilde{g}}_{\mathbb{H}^2(\mathcal{S}^\ast)} \lesssim \norm{g}_{\mathbb{H}^2(\mathcal{S}^\ast)}.
    \end{align*}
    Moreover, as $0$ belongs to the closure of the image of $\widehat{\mathscr{P}}_{\pi-\theta}^2$ (since $|\widehat{\mathscr{P}}_{\pi-\theta}(x+\ii y)| \ra 0$ as $|x| \ra \infty$), we have $\lambda \neq 0$ and therefore the function $\tilde{f}$ defined in~\eqref{eq:solution}  satisfies 
\begin{align*}
    \norm{\tilde{f}}_{\mathbb{H}^2(\mathcal{S}^\ast)} \lesssim \norm{g}_{\mathbb{H}^2(\mathcal{S}^\ast)} + \norm{f}_{\mathbb{H}^2(\mathcal{S}^\ast)}.
\end{align*}
Consequently, the inverse map defined in~\eqref{eq:solution} is continuous and $\lambda \not \in \sigma(T_{21})$, which concludes the proof of the first statement.

To show that any $\lambda \in \{ \overline{\lambda} \in \C : \lambda \in \mathscr{P}^2_{\pi-\theta}(\mathcal{S}^\ast)\}$ is an eigenfunction of $T_{21}^\ast$, we can use Lemma~\ref{lem:Poissonproperties}~\ref{it:reproducingkernel}. Indeed, let $\lambda = \widehat{\mathscr{P}}_{\pi-\theta}(z_0)^2$ for some $z_0 \in \mathcal{S}^\ast$, then, from~\eqref{eq:solution}, we see that
\begin{align*}
    \inner{(0,\mathscr{T}_{z_0} \widehat{\mathscr{P}}_{\sigma+\pi}), (\lambda- T_{21})(\tilde{f},\tilde{g})}_{\mathbb{H}^2_{\sigma}\times \mathbb{H}^2_{\sigma+\pi}} &= \inner{\mathscr{T}_{z_0} \widehat{\mathscr{P}}_{\sigma+\pi}, (\lambda - \widehat{\mathscr{P}}_{\pi-\theta}^2) \tilde{g}}_{\sigma+\pi} \\
    &= \frac{\sin(\sigma+\pi)}{2\pi} \left(\lambda - \widehat{\mathscr{P}}_{\pi-\theta}(z_0)^2\right) \tilde{g}(z_0) = 0.
\end{align*}
As the above holds for any $\tilde{f},\tilde{g} \in \mathbb{H}^2(\mathcal{S}^\ast)$, the vector $(0,\mathscr{T}_{z_0} \widehat{\mathscr{P}}_{\sigma+\pi}) \in (\mathbb{H}^2(\mathcal{S}^\ast))^2$ is an eigenfunction of the adjoint of $T_{21}$ with eigenvalue $\overline{\lambda} = \overline{\widehat{\mathscr{P}}_{\pi-\theta}(z_0)^2}$. 

To obtain the expression in~\eqref{eq:expliciteigenfunction}, we can assume w.l.o.g. that  $\Omega_1$ and $\Omega_2$ are in the configuration from Section~\ref{sec:bipolartransform}. Hence, by Theorem~\ref{thm:DtDHardy},
\begin{align*}
   \left(0, \lambda \mathscr{T}_{z_0} \widehat{\mathscr{P}}_{\sigma+\pi}\right) = T_{21}^\ast(0, \mathcal{T}_{z_0} \widehat{\mathscr{P}}_{\sigma+\pi}) = \left(0,\mathcal{F} W_{\sigma+\pi}^2 \Psi_{\sigma+\pi}^\# \gamma_1 \gamma_2 (\Psi_{\sigma+\pi}^\#)^{-1}W_{\sigma+\pi}^{-2} \mathcal{F}^{-1} \mathscr{T}_{z_0} \widehat{\mathscr{P}}_{\sigma+\pi}\right)
\end{align*}
 where $W_{\sigma+\pi}$ is the weight multiplication defined in~\eqref{eq:weight} and $\Psi^\#_{\sigma+\pi}$ is the pullback with respect to the bipolar transform defined in~\eqref{eq:pullback}. Consequently, the function
\begin{align*}
    f_{z_0}(\tau) = W_{\sigma+\pi}^{-2} \mathcal{F}^{-1} \mathcal{T}_{z_0} \widehat{\mathscr{P}}_{\sigma+\pi} = \ee^{\ii \overline{z_0} \tau} 
\end{align*}
is an eigenfunction of the $\Psi^\# \gamma_1\gamma_2 (\Psi^{\#})^{-1}$. To complete the proof, we can now compute the pushforward of $f_{z_0}(\tau)$ with respect to the bipolar transform defined in~\eqref{eq:bipolardef}.
\end{proof}

\subsection{Proof of Theorem~\ref{thm:DtDests}~\ref{it:ext and  int}} Next, we prove a improved version of the estimates in Theorem~\ref{thm:DtDests}~\ref{it:ext and  int}. These improved estimates will be necessary for the proof of Theorem~\ref{thm:DtDests}~\ref{it:numericalrange}.

\begin{theorem}[Interior and exterior estimates] \label{thm:interior_exterior} For any $\theta \in (0,\pi)$, let $P_\theta : \mathbb{H}^2(\mathcal{S}^\ast) \rightarrow \mathbb{H}^2(\mathcal{S}^\ast)$ denote the multiplication operator
\begin{align*}
(P_\theta h)(z) = \widehat{\mathscr{P}}_\theta(z) h(z).
\end{align*}
Then for any $\sigma \in  (0,\pi)$ with $\sigma+\theta \leq \pi$, there exists $h^{\rm int} ,h^{\rm ext} \in \mathbb{H}^2(\mathcal{S}^\ast)$ (depending on $\theta,\sigma$) such that 
for any $h\in \mathbb{H}^2(\mathcal{S}^\ast)$ we have
\begin{align}
    &\norm{P_\theta h}_{\mathbb{H}^2_{\sigma}}^2 \leq\cos(\theta/2)^2  \norm{ P_{h^{\rm ext}}^\perp h}_{\mathbb{H}^2_{\sigma+\theta}}^2 + |\inner{h,h^{\rm ext}}_{\mathbb{H}^2_{\sigma+\theta}}|^2,\label{eq:onesidedest1} \\
    &\norm{P_{\pi-\theta} h}_{\mathbb{H}^2_{\sigma+\pi}}^2 \leq \sin(\theta/2)^2 \norm{P^\perp_{h^{\rm int}}h}_{\mathbb{H}^2_{\sigma+\theta}}^2 + |\inner{h,h^{\rm int}}_{\mathbb{H}^2_{\sigma+\theta}}|^2, \label{eq:onesidedest2}
\end{align}
where $P^\perp_{h^{\rm int/ext}}$ denotes the orthogonal projection on the orthogonal complement of $h^{\rm int/ext}$. Moreover, we have the bounds
\begin{align}
0 \leq &\norm{P_\theta}_{\mathbb{H}^2_{\sigma+\theta}\rightarrow \mathbb{H}^2_{\sigma}}^2 -\cos(\theta/2)^2\leq \left(\cos(\sigma+\theta)\right)_+ \frac{\sin(\theta)}{\sin(\sigma+\theta)} \frac{\pi-\theta}{\pi} \frac{\pi-\sigma-\theta}{\pi} \label{eq:onesidedupper1} \\
0 \leq &\norm{P_{\pi-\theta}}_{\mathbb{H}^2_{\sigma+\theta} \rightarrow \mathbb{H}^2_{\sigma+\pi}}^2 -\sin(\theta/2)^2 \leq  \left(\frac{\cos(\sigma+\theta) + \cos(\sigma) \frac{\theta^2}{\pi^2 \sin(\theta/2)^2}}{1+\cos(\theta) \frac{\theta^2}{\pi^2\sin(\theta/2)^2}}\right)_+ \frac{\sin(\theta)}{\sin(\theta+\sigma)} \frac{\theta}{\pi}\frac{\pi-\theta-\sigma}{\pi}, \label{eq:onesidedupper2} 
\end{align}
where $a_+ = \max \{a, 0\}$.
\end{theorem}

\begin{proof} The proof here is similar to the proof of Theorem~\ref{thm:generalbound}. Hence, for the sake of brevity, we give the full details for the proof of~\eqref{eq:onesidedest2} and~\eqref{eq:onesidedupper2}, and only briefly sketch the proof of~\eqref{eq:onesidedest1} and~\eqref{eq:onesidedupper1}. First, we note that
\begin{align}
    \norm{h}_{\sigma+\theta}^2 - \frac{1}{\sin(\theta/2)^2} \norm{\widehat{\mathscr{P}}_{\pi-\theta} h}_{\sigma+\pi}^2 =& \frac{1}{2\pi} \int_\R \left(1-\frac{|\widehat{\mathscr{P}}_{\pi-\theta}(x+\ii/2)|^2}{\sin(\theta/2)^2}\right) \frac{|h_{1/2}(x)|^2 + |h_{-1/2}(x)|^2}{2} \mathrm{d} x \nonumber \\
    &- \frac{1}{2\pi} \int_\R \left(\cos(\sigma) \frac{\widehat{\mathscr{P}}_{\pi-\theta}(x)^2}{\sin(\theta/2)^2} + \cos(\sigma+\theta)\right) |h(x)|^2 \mathrm{d} x  \label{eq:someequation2}
 \end{align}
 As we did in the proof of Theorem~\ref{thm:generalbound}, we can now use the pointwise identity~\eqref{eq:pointwiseid} and apply the residue theorem in~\eqref{eq:residuethm} with 
 \begin{align*}
     F(z) \coloneqq  1-\frac{\widehat{\mathscr{P}}_{\pi-\theta}(z) \widehat{\mathscr{P}}_{\pi-\theta}(z-\ii)}{\sin(\theta/2)^2}, 
\end{align*}
which satisfies
\begin{align*}
     \mathrm{Re} \, F(x) = 1+\cos(\theta) \frac{\widehat{\mathscr{P}}_{\pi-\theta}(x)^2}{\sin(\theta/2)^2}  \quad \mbox{and}\quad \mathrm{Res}(F,0) = - \frac{\ii}{\pi} \frac{\theta}{\pi} \frac{\sin(\theta)}{\sin(\theta/2)^2}, 
 \end{align*}
 to re-write the right-hand side of~\eqref{eq:someequation2} as
 \begin{align*}
     \frac{1-\alpha}{2\pi} \int_\R \left(1-\frac{|\widehat{\mathscr{P}}_{\pi-\theta}(x+\ii/2)|^2}{\sin(\theta/2)^2}\right) \frac{|h_{1/2}(x)|^2 + |h_{-1/2}(x)|^2}{2} \mathrm{d} x\\
     + \frac{\alpha}{2\pi} \int_\R \left(1-\frac{|\widehat{\mathscr{P}}_{\pi-\theta}(x+\ii/2)|^2}{\sin(\theta/2)^2}\right) \frac{|h_{1/2}(x)-h_{-1/2}(x)|^2}{2} \mathrm{d} x \\
     + \frac{1}{2\pi} \int_\R\left(\alpha-\cos(\sigma+\theta) + \left(\alpha\cos(\theta) - \cos(\sigma)\right) \frac{|\widehat{\mathscr{P}}_{\pi-\theta}(x)|^2}{\sin(\theta/2)^2}\right) |h_0(x)|^2 \mathrm{d} x\\
     - \frac{\alpha}{2\pi} \frac{\sin(\theta)}{\sin(\theta/2)^2} \frac{\theta}{\pi} |h(0)|^2,
 \end{align*}
 for any $\alpha \in \R$. Since $\widehat{\mathscr{P}}_{\pi-\theta}(x) \leq \sin(\theta/2)$ by Lemma~\ref{lem:Poissonproperties}~\ref{it:Poissonmaxima}, the first two terms are non-negative as long as $0\leq \alpha \leq 1$. To have a non-negative integrand also for the third integral, we need
 \begin{align}
     \alpha -\frac{\sin(\theta/2)^2\cos(\sigma+\theta) + \cos(\sigma) \widehat{\mathscr{P}}_{\pi-\theta}(x)^2}{\sin(\theta/2)^2 + \cos(\theta) \widehat{\mathscr{P}}_{\pi-\theta}(x)^2} \geq 0 \quad \mbox{for any $x\in \R$.} \label{eq:differenceeq}
 \end{align}
 
 We now claim that the smallest possible $\alpha \in \R$ for~\eqref{eq:differenceeq} to hold is given by
 \begin{align}
     \alpha \coloneqq \frac{\sin(\theta/2)^2\cos(\sigma+\theta) + \cos(\sigma) \widehat{\mathscr{P}}_{\pi-\theta}(0)^2}{\sin(\theta/2)^2 + \cos(\theta) \widehat{\mathscr{P}}_{\pi-\theta}(0)^2} = \frac{\cos(\sigma+\theta) + \cos(\sigma) \frac{\theta^2}{\pi^2 \sin(\theta/2)^2}}{1+\cos(\theta) \frac{\theta^2}{\pi^2 \sin(\theta/2)^2}}. \label{eq:alphadef}
 \end{align}
 To see this, note that the difference in~\eqref{eq:differenceeq} with $\alpha$ given by~\eqref{eq:alphadef} is
 \begin{align*}
     &\sin(\theta/2)^2 \frac{\left( \cos(\sigma)-\cos(\sigma+\theta)\cos(\theta)\right)\left(\widehat{\mathscr{P}}_{\pi-\theta}(0)^2 - \widehat{\mathscr{P}}_{\pi-\theta}(x)^2\right)}{\left(\sin(\theta/2)^2 + \cos(\theta) \widehat{\mathscr{P}}_{\pi-\theta}(0)^2\right)\left(\sin(\theta/2)^2 + \cos(\theta) \widehat{\mathscr{P}}_{\pi-\theta}(x)^2\right)} \\
     &= \sin(\theta/2)^2 \frac{\sin(\sigma+\theta) \sin(\theta)\left(\widehat{\mathscr{P}}_{\pi-\theta}(0)^2 - \widehat{\mathscr{P}}_{\pi-\theta}(x)^2\right)}{\left(\sin(\theta/2)^2 + \cos(\theta) \widehat{\mathscr{P}}_{\pi-\theta}(0)^2\right)\left(\sin(\theta/2)^2 + \cos(\theta) \widehat{\mathscr{P}}_{\pi-\theta}(x)^2\right)} .
 \end{align*}
 Since $\sin(\sigma+\theta) \geq 0$ (as $0 \leq \sigma+\theta \leq \pi$) and $\widehat{\mathscr{P}}_{\pi-\theta}(x)^2 \leq \widehat{\mathscr{P}}_{\pi-\theta}(0)^2 = \theta^2/\pi^2  \leq \sin(\theta/2)^2$, the above expression is exactly $0$ if $x=0$ and positive otherwise. This proves our claim.
 
 We have thus shown that
 \begin{align*}
     \norm{h}_{\sigma+\theta}^2 - \frac{1}{\sin(\theta/2)^2} \norm{\widehat{\mathscr{P}}_{\pi-\theta} h}_{\sigma+\pi}^2 \geq \alpha_+ \frac{\sin(\theta)}{\sin(\theta/2)^2} \frac{\theta}{2\pi^2} |h(0)|^2 
 \end{align*}
 with $\alpha$ defined by~\eqref{eq:alphadef}. If $\alpha$ is non-positive, this proves~\eqref{eq:onesidedupper2} as well as~\eqref{eq:onesidedest2} with $h^{\rm ext} =0$. If $\alpha >0$, then we note that the subspace of functions with $h(0) \neq 0$ has codimension $1$. Thus by the min-max principle, the part of the spectrum of the symmetric operator $(P_\theta)^\ast P_\theta$ inside the interval $(\sin(\theta/2)^2, \infty)$ contains at most one simple eigenvalue. This implies the existence of at most one non-trivial $h^{\rm ext}$ such that~\eqref{eq:onesidedest1} holds. To prove the upper bound in~\eqref{eq:onesidedupper2} we can use Lemma~\ref{lem:Poissonproperties}~\ref{it:reproducingkernel} and the Cauchy-Schwarz inequality, as we did in~\eqref{eq:CS}.

For~\eqref{eq:onesidedest2} and~\eqref{eq:onesidedupper2}, we can follow the same steps to show that
\begin{align*}
\norm{h}_{\sigma+\theta}^2 - \frac{\norm{\widehat{\mathscr{P}}_{\theta} h}_{\sigma}^2}{\cos(\theta/2)^2}  = \frac{1-\alpha}{2\pi} \int_\R \left(1-\frac{|\widehat{\mathscr{P}}_{\theta}(x+\ii/2)|^2}{\cos(\theta/2)^2}\right) \frac{|h_{1/2}(x)|^2 + |h_{-1/2}(x)|^2}{2} \mathrm{d} x\\
     + \frac{\alpha}{2\pi} \int_\R \left(1-\frac{|\widehat{\mathscr{P}}_{\theta}(x+\ii/2)|^2}{\cos(\theta/2)^2}\right) \frac{|h_{1/2}(x)-h_{-1/2}(x)|^2}{2} \mathrm{d} x \\
     + \frac{1}{2\pi} \int_\R\left(\alpha-\cos(\sigma+\theta) - \left(\alpha\cos(\theta) - \cos(\sigma)\right) \frac{|\widehat{\mathscr{P}}_{\theta}(x)|^2}{\cos(\theta/2)^2}\right) |h_0(x)|^2 \mathrm{d} x\\
     - \frac{\alpha}{2\pi} \frac{\sin(\theta)}{\cos(\theta/2)^2} \frac{\pi-\theta}{\pi} |h(0)|^2.
\end{align*}
As before, this already implies~\eqref{eq:onesidedest1}. Moreover, one can show that the minimum $\alpha \in [0,1]$ for which the integrand in the third integral is non-negative is $\alpha = \cos(\sigma+\theta)_+$. Hence, the same argument from before combining Cauchy-Schwarz with Lemma~\ref{lem:Poissonproperties}~\ref{it:reproducingkernel} (see~\eqref{eq:CS}) yields~\eqref{eq:onesidedupper2}. 

The lower bounds in~\eqref{eq:onesidedupper1} and~\eqref{eq:onesidedupper2} can be shown again by using the localization procedure in Lemma~\ref{lem:localizationlemma}. As the same arguments were already used in the proof of Theorem~\ref{thm:generalbound} and will be used again in the proof of Lemma~\ref{lem:coercivity} below, we skip the details here.
\end{proof}

\subsection{Proof of Theorem~\ref{thm:DtDests}~\ref{it:numericalrange}} We are now ready to prove Theorem~\ref{thm:DtDests}~\ref{it:numericalrange}. By Theorem~\ref{thm:DtDHardy}, it suffices to bound (from below and from above) 
\begin{align*}
    \sup_{\substack{(h,k) \in (\mathbb{H}^2)^2\\ (u,v) \in (\mathbb{H}^2)^2}} \frac{|\inner{(u,v), T_1(h,k)}_{\mathbb{H}^2_{\sigma+\theta}\times \mathbb{H}^2_{\sigma+\theta+\pi}}+\inner{(h,k), T_2(u,v)}_{\mathbb{H}^2_\sigma \times \mathbb{H}^2_{\sigma+\pi}}|}{\norm{h}_{\sigma}^2 + \norm{k}_{\sigma+\pi}^2 + \norm{u}_{\sigma+\theta}^2 + \norm{v}_{\sigma+\theta+\pi}^2},
\end{align*}
where $T_1$ and $T_2$ are the operators defined in~\eqref{eq:T1def} and \eqref{eq:T2def}. Hence, the result follows from the following lemma. 

\begin{lemma}[Numerical radius estimate] \label{lem:coercivity} Let $\theta,\sigma \in (0,\pi)$ with $\theta+\sigma \leq \pi$. Then
\begin{align}
    \frac{1+\sin(\theta/2)}{2} \leq \!\! \sup_{\substack{(h,k) \in (\mathbb{H}^2)^2\\ (u,v) \in (\mathbb{H}^2)^2}} \frac{|\inner{(u,v), T_1(h,k)}_{\mathbb{H}^2_{\sigma+\theta}\times \mathbb{H}^2_{\sigma+\theta+\pi}}+\inner{(h,k), T_2(u,v)}_{\mathbb{H}^2_\sigma \times \mathbb{H}^2_{\sigma+\pi}}|}{\norm{h}_{\sigma}^2 + \norm{k}_{\sigma+\pi}^2 + \norm{u}_{\sigma+\theta}^2 + \norm{v}_{\sigma+\theta+\pi}^2} \leq f(\theta), \label{eq:coercivitybounds}
\end{align}
where
\begin{align*}
    f(\theta) = \frac{1}{2} \left(1+\sin(\theta/2)^2 + g(\theta) + \sqrt{4\sin(\theta/2)^2 + h(\theta)}\right)^{\frac{1}{2}},
\end{align*}
with
\begin{align}
&g(\theta) \coloneqq \left(\frac{\alpha(\theta)_+}{2} \frac{\theta}{\pi} +  \frac{\cos(\theta)_+}{2}\frac{\pi-\theta}{\pi}\right)\frac{\pi-\theta}{\pi} -\sin(\theta/2)^2\left(1 - \sqrt{1 + \frac{\alpha(\theta)_+}{\sin(\theta/2)^2} \frac{\theta}{\pi} \frac{\pi-\theta}{\pi}}\right),\label{eq:gfunction} \\
        &\alpha(\theta) = \frac{ \pi^2 \sin(\theta/2)^2 \cos(\theta) + \theta^2}{\pi^2 \sin(\theta/2)^2 + \theta^2 \cos(\theta)}, \nonumber \\
	&h(\theta) \coloneqq g(\theta)^2 +  2g(\theta) \left(1+\sin(\theta/2)^2\right) - \cos(\theta/2)^2 \cos(\theta)_+ \frac{(\pi-\theta)^2}{\pi^2} , \nonumber
\end{align}
and $a_+ = \max (0,a)$.
\end{lemma}

\begin{proof}
    Let us start by proving the lower bound. For this, we use again the localization Lemma~\ref{lem:localizationlemma}. So let $f\in \mL^2(\R)$ and define the function
\begin{align*}
    f^\delta(z) \coloneqq \left(m D_\delta f\right)(z), 
\end{align*}
where $m$ is given by~\eqref{eq:moperator}. Then we choose as trial states in~\eqref{eq:coercivitybounds} the functions
\begin{align*}
(h,k) \coloneqq  \left(c \widehat{\mathscr{P}}_\theta f^\delta, f^\delta\right)
\quad \mbox{and}\quad (u,v) \coloneqq \left( f^\delta, c \widehat{\mathscr{P}}_\theta f^\delta\right) \quad\mbox{with} \quad c \coloneqq \frac{1}{\sin(\theta/2)+1}.
\end{align*}
From the symmetry 
\begin{align*}
    \mathrm{Re} \, \widehat{\mathscr{P}}_{\pi-\theta}(x+\ii/2) = \mathrm{Re}\, \widehat{\mathscr{P}}_{\pi-\theta}(x-\ii/2), 
\end{align*}
we see that the real part of the numerator in~\eqref{eq:coercivitybounds} for this choice of trial states is equal to
\begin{align*}
    \mathrm{Re} &\left\{\inner{u, \widehat{\mathscr{P}}_{\pi-\theta} k}_{\sigma+\theta} +\inner{k,\widehat{\mathscr{P}}_{\pi-\theta} u}_{\sigma+\pi} + \inner{v,\widehat{\mathscr{P}}_\theta  k}_{\sigma+\theta+\pi} + \inner{h,\widehat{\mathscr{P}}_\theta u}_{\sigma}\right\} = \\
    =& \frac{1}{2\pi} \int_{\R} 2\left( \mathrm{Re} \left\{ \widehat{\mathscr{P}}_{\pi-\theta}(x+\ii/2)\right\} +  c |\widehat{\mathscr{P}}_\theta(x+\ii/2)|^2 \right)\frac{|f^\delta_{1/2}(x)|^2 +|f^\delta_{-1/2}(x)|^2}{2}\mathrm{d}x \\
    &-\frac{1}{2\pi} \int_\R \left(\cos(\sigma+\theta) - \cos(\sigma)\right) \left(\widehat{\mathscr{P}}_{\pi-\theta}(x) - c \widehat{\mathscr{P}}_{\theta}(x)^2\right) |f^\delta(x)|^2 \mathrm{d} x.
\end{align*}
Thus by re-scaling in $x$ by $\delta$ and applying~\eqref{eq:secondtransformidentity}, the second integral vanishes in the limit $\delta\downarrow 0$. Moreover, by~\eqref{eq:firsttransformidentity},~\eqref{eq:secondtransformidentity}, and dominated convergence, the first integral converge towards
\begin{align}
    \frac{1}{\pi} \left( \mathrm{Re} \widehat{\mathscr{P}}_{\pi-\theta}(\ii/2) + c|\widehat{\mathscr{P}}_\theta(\ii/2)|^2\right)\left(\norm{f}_{\mL^2(\R)}^2 + \norm{Hf}_{\mL^2(\R)}^2\right) &= \frac{2}{\pi} \left(\sin(\theta/2) + \frac{\cos(\theta/2)^2}{1+\sin(\theta/2)}\right) \norm{f}_{\mL^2(\R)}^2 \nonumber \\
    &= \frac{2}{\pi} \norm{f}_{\mL^2(\R)}^2. \label{eq:converged}
\end{align}
By looking now to the denominator in~\eqref{eq:coercivitybounds} and using Lemma~\ref{lem:localizationlemma}, we find that
\begin{align}
    \lim_{\delta \downarrow 0} \norm{h}_{\sigma}^2 + \norm{k}_{\sigma+\pi}^2+\norm{u}_{\sigma+\theta}^2 + \norm{v}_{\sigma+\theta+\pi}^2 &= \lim_{\delta \downarrow 0} c^2 \norm{\widehat{\mathscr{P}}_\theta f^\delta}_{\sigma}^2 + \norm{f^\delta}_{\sigma+\pi}^2 + \norm{f^\delta}_{\sigma+\theta}^2 + c^2\norm{\widehat{\mathscr{P}}_\theta f^\delta}_{\sigma+\theta+\pi}^2 \nonumber \\
    &= \frac{1}{\pi} \left( c^2 \widehat{\mathscr{P}}_\theta(\ii/2)^2 + 1 \right)\left(\norm{f}_{\mL^2(\R)}^2 + \norm{Hf}_{\mL^2(\R)}^2\right)\nonumber \\
    &= \frac{2}{\pi}\left(\frac{\cos(\theta/2)^2}{\left(1+\sin(\theta/2)\right)^2} + 1\right) \norm{f}_{\mL^2(\R)}^2 =  \frac{4c}{\pi} \norm{f}_{\mL^2(\R)}^2. \label{eq:denominator}
\end{align}
As the real part gives a lower bound on the absolute value, estimate~\eqref{eq:denominator} together with~\eqref{eq:converged} yields the lower bound in~\eqref{eq:coercivitybounds}.

To prove the upper bound, we first note that by Cauchy-Schwarz and the trivial bound $\norm{P_\theta f} \leq \norm{P_\theta} \norm{f}$ we have
\begin{multline}
    |\inner{u, \widehat{\mathscr{P}}_{\pi-\theta} k}_{\sigma+\theta} +\inner{k,\widehat{\mathscr{P}}_{\pi-\theta} u}_{\sigma+\pi} + \inner{v,\widehat{\mathscr{P}}_\theta  k}_{\sigma+\theta+\pi} + \inner{h,\widehat{\mathscr{P}}_\theta u}_{\sigma}|  \\
\leq c_1 \norm{u}_{\sigma+\theta} \norm{k}_{\sigma+\pi}+ c_2\norm{u}_{\sigma+\theta} \norm{h}_{\sigma} +  c_3\norm{k}_{\sigma+\pi} \norm{v} , \label{eq:somemiddleest}
\end{multline}
for any positive constants $c_1, c_2, c_3 \geq 0$ satisfying
\begin{align*}
    c_1 \geq \norm{P_{\pi-\theta}}_{\mathbb{H}^2_{\sigma+\theta} \rightarrow \mathbb{H}^2_{\sigma+\pi}} + \norm{P_{\pi-\theta}}_{\mathbb{H}^2_{\sigma+\pi} \rightarrow \mathbb{H}^2_{\sigma+\theta}},\quad c_2 \geq \norm{P_\theta}_{\mathbb{H}^2_{\sigma+\theta} \rightarrow \mathbb{H}^2_{\sigma}}, \quad\mbox{and}\quad c_3 \geq \norm{P_\theta}_{\mathbb{H}^2_{\sigma+\pi} \rightarrow \mathbb{H}^2_{\sigma+\theta+\pi}}. 
\end{align*}
From Young's inequality, the right-hand side of~\eqref{eq:somemiddleest} can be bounded by
\begin{align*}
    \left(\frac{c_1}{2} \frac{1}{\epsilon} + \frac{c_2^2}{4 \delta}\right) \norm{u}_{\sigma+\theta}^2 + \left(\frac{c_1}{2}\epsilon + \frac{c_3^2}{4 \delta}\right) \norm{k}_{\sigma+\pi}^2 + \delta \left(\norm{v}_{\sigma+\theta+\pi}^2 + \norm{h}_\sigma^2\right),
\end{align*}
for any $\delta, \epsilon>0$. Balancing the terms in the above expression, i.e., solving for $\delta>0$ and $\epsilon>0$ the equations
\begin{align*}
    \frac{c_1}{2} \frac{1}{\epsilon} + \frac{c_3^2}{4 \delta} = \delta = \frac{c_1}{2} \epsilon + \frac{c_2^2}{4 \delta},
\end{align*}
we find that $\delta$ is given by
\begin{align}
    \delta(c_1,c_2,c_3) = \left(\frac{c_1^2 + c_2^2 + c_3^2 + \sqrt{(c_1^2+c_2^2+c_3^2)^2 - 4 c_2^2 c_3^2}}{8} \right)^{\frac12}. \label{eq:deltaformula}
\end{align}
Next, we note that, since $\cos(2\pi - \phi) = \cos(\phi)$ for any $\phi \in \R$, we have
\begin{align*}
\norm{P_{\pi-\theta}}_{\mathbb{H}^2_{\sigma+\pi} \rightarrow \mathbb{H}^2_{\sigma+\theta}} = \norm{P_{\pi-\theta}}_{\mathbb{H}^2_{\pi-\sigma} \rightarrow \mathbb{H}^2_{2\pi-\sigma-\theta}} \quad \mbox{and}\quad \norm{P_\theta}_{\mathbb{H}^2_{\sigma+\pi} \rightarrow \mathbb{H}^2_{\sigma+\theta+\pi}}  = \norm{P_\theta}_{\mathbb{H}^2_{\pi-\sigma} \rightarrow \mathbb{H}^2_{\pi -\sigma-\theta}}
\end{align*}
Therefore, by the upper bounds in~\eqref{eq:onesidedupper1} and~\eqref{eq:onesidedupper2}, we can choose $c_1$, $c_2$ and $c_3$ as
\begin{align*}
&c_1(\theta,\sigma) = b(\theta,\sigma) + b(\theta,\pi-\theta-\sigma), \\
&c_2(\theta,\sigma) =  \left(\cos(\theta/2)^2 + \left( \frac{\sin(\theta)}{\tan(\sigma+\theta)}\right)_+ \frac{\pi-\theta}{\pi}\frac{\pi-\theta-\sigma}{\pi}\right)^{\frac12}, \\
&c_3(\theta,\sigma) = c_2(\theta,\pi-\theta-\sigma), 
\end{align*}
where
\begin{align*}
    b(\theta,\sigma) \coloneqq \left(\sin(\theta/2)^2 + \left(\frac{\sin(\theta/2)^2 \cos(\sigma+\theta) + \cos(\sigma) \frac{\theta^2}{\pi^2}}{\sin(\theta/2)^2+\cos(\theta) \frac{\theta^2}{\pi^2}}\right)_+ \frac{\sin(\theta)}{\sin(\sigma+\theta)} \frac{\theta}{\pi} \frac{\pi-\theta-\sigma}{\pi} \right)^{\frac12}.
\end{align*}
If we now use these expressions back in~\eqref{eq:deltaformula} we obtain a $(\theta,\sigma)$-dependent upper bound for the left-hand side of~\eqref{eq:coercivitybounds}. To obtain the $\sigma$-independent upper bound stated there, we fix $\theta \in (0,\pi)$ and note that, from the symmetric relations 
\begin{align*}
    \delta(c_1,c_2,c_3) = \delta(c_1,c_3,c_2), \quad c_1(\theta,\sigma) = c_1(\theta,\pi-\theta-\sigma)\quad \mbox{and}\quad c_2(\theta,\sigma) = c_3(\theta,\pi-\theta-\sigma),
\end{align*}
it suffices to consider $\sigma$ in the interval $(0,\frac{\pi-\theta}{2})$. For any $\sigma$ in this interval, we note that
\begin{align*}
    c_3(\theta,\sigma) = \cos(\theta/2) \quad \mbox{and} \quad b(\theta,\pi-\theta-\sigma) = \sin(\theta/2). 
\end{align*}
Moreover, the functions $c_2(\theta,\sigma)$ and $c_1(\theta,\sigma)$ defined above are monotonically decreasing for $\sigma \in [0,\frac{\pi-\theta}{2}]$. Indeed, for $c_2(\theta,\sigma)$, this follows from the fact that $(\pi-\theta-\sigma)/\tan(\sigma+\theta)$ is decreasing for $\sigma+\theta \in [0,\pi]$. For $c_1(\theta,\sigma)$, this follows from the fact that
\begin{align*}
    \frac{\sin(\theta/2)^2 \cos(\sigma+\theta) + \cos(\sigma)\frac{\theta^2}{\pi^2}}{\sin(\theta/2)^2 + \cos(\theta) \frac{\theta^2}{\pi^2}} \frac{\pi-\theta-\sigma}{\sin(\sigma+\theta)} = \frac{\pi-\theta-\sigma}{\tan(\sigma+\theta)} + \frac{(\pi-\theta-\sigma)\theta^2}{\pi^2 \sin(\theta/2)^2 + \cos(\theta) \theta^2}
\end{align*}
 is also decreasing for $\sigma+\theta \in (0,\pi)$. Hence, the maximum of $c_1(\theta,\sigma)$, $c_2(\theta,\sigma)$, and $c_3(\theta,\sigma)$ for $\sigma \in [0,\frac{\pi-\theta}{2})$ is simultaneously achieved at $\sigma = 0$. Consequently,
\begin{align}
    \delta(\theta) \coloneqq &\sup_{\sigma \in (0,\frac{\pi-\theta}{2})} \delta\left(c_1(\theta,\sigma),c_2(\theta,\sigma),c_3(\theta,\sigma)\right) = \delta\left(c_1(\theta,0),c_2(\theta,0),c_3(\theta,0)\right). \label{eq:maxdelta}
\end{align}
To conclude, note that
\begin{align*}
    c_1(\theta,0)^2 + c_2(\theta,0)^2 + c_3(\theta,0)^2 &= 4\sin(\theta/2)^2 + 2 \cos(\theta/2)^2 + 2g(\theta)= 2\left(1+\sin(\theta/2)^2 + g(\theta)\right),
\end{align*}
where $g(\theta)$ is defined in~\eqref{eq:gfunction}, and
\begin{align*} c_2(\theta,0)^2c_3(\theta.0)^2 = \cos(\theta/2)^4 + \cos(\theta/2)^2\cos(\theta)_+ \frac{(\pi-\theta)^2}{\pi^2},
\end{align*}
and therefore the result follows from~\eqref{eq:maxdelta} and the formula in~\eqref{eq:deltaformula}.
\end{proof}

\addtocontents{toc}{\protect\setcounter{tocdepth}{-1}}
\section*{Acknowledgements}
The authors thank Xavier Claeys, Yvon Maday, and Zhuoyao Zeng for helpful discussions. AJ and BS acknowledge support from the \emph{Deutsche Forschungsgemeinschaft} (DFG, German Research Foundation) - Project number 440641818. BS and TC acknowledge funding by the \emph{Deutsche Forschungsgemeinschaft} (DFG, German Research Foundation) - Project number 442047500 through the Collaborative Research Center "Sparsity and Singular Structures" (SFB 1481).

\addtocontents{toc}{\protect\setcounter{tocdepth}{2}}
\appendix
\section{Proof of auxiliary lemmas}
\label{app:Poisson}
In this section we present the proof of some of the technical lemmas that were used throughout Section~\ref{sec:proof}. We start with a proof of Lemma~\ref{lem:Poissonproperties}.

\begin{proof}[Proof of Lemma~\ref{lem:Poissonproperties}]

To prove~\ref{it:Poissonmaxima} we note that, since $\mathscr{P}_\theta(x)$ is non-negative and even, from the Fourier transform formula we have
\begin{align*}
|\widehat{\mathscr{P}}_\theta(x+\ii y)| \leq \int_\R \ee^{y k} \mathscr{P}_\theta(k) \mathrm{d} k = \int_\R \mathscr{P}_\theta(k) \cosh(yk) \mathrm{d}k = \widehat{\mathscr{P}}_\theta(\ii y)\quad \mbox{for any $x\in \R$,} 
\end{align*}
with strict inequality if $x\neq 0$. Since $\cosh(yk) <\cosh(y'k)$ for any $|y'k| > |y k|$, we have $\widehat{\mathscr{P}}_\theta(\ii y') > \widehat{\mathscr{P}}_\theta(\ii y)$ in this case. The result then follows by taking the limit $|y| \uparrow 1/2$.

To prove item~\ref{it:combinedsupremumnorm}, we first note that the rightmost inequality is elementary. Indeed, it follows from the fact that $f(\theta) = 1-2\theta/\pi-\cos(\theta)$ has only a local minimum in the interval $(0,\pi/2)$ and therefore the maximum is achieved at $f(0) = f(\pi/2) = 0$. The leftmost inequality in~\ref{it:combinedsupremumnorm} follows from the formula for the Poisson kernel~\eqref{eq:Poissondef} and the simple inequality $|\sinh(x)| \leq |\sinh( x')|$ for $|x| \leq |x'|$. Hence, it remains to show that
\begin{align*}
    \widehat{\mathscr{P}}_\theta(x)^2 - \widehat{\mathscr{P}}_{\pi-\theta}(x)^2 = \frac{\cosh(2(\pi-\theta)x)-\cosh(2\theta x)}{\cosh(2\pi x)-1} \leq 1-\frac{2\theta}{\pi},
\end{align*}
where we used the identity $2\sinh(a)^2 = \cosh(2a)^2 - 1$ for any $a\in \R$. To prove this inequality, it suffices to show that the function
\begin{align*}
    f_\theta(x) \coloneqq \cosh((\pi-\theta)x) - \cosh(\theta x) - \left(1-\frac{2\theta}{\pi}\right) \left(\cosh(\pi x)-1\right)
\end{align*}
is non-positive for $x\geq 0$. For this, note that
\begin{align}
   f_\theta^{(n)}(0) = \begin{dcases} (\pi-\theta)^n - \theta^n - \left(1-\frac{2\theta}{\pi}\right) \pi^n, \quad &\mbox{for $n \in \N$ even,} \\
   0 &\mbox{ for $n\in \N$ odd.} \end{dcases} \label{eq:derivative of f}
\end{align}
Since $\theta\leq \pi-\theta$ (as $\theta\leq \pi/2$), for $n\geq 2$ we have
\begin{align*}
    (\pi-\theta)^n - \theta^n - \left(1-\frac{2\theta}{\pi}\right) \pi^n &=\pi \left((\pi-\theta)^{n-1} - \theta^{n-1} -\left(1-\frac{2\theta}{\pi}\right) \pi^{n-1}\right) - \theta (\pi-\theta)^{n-1} + \theta^{n-1}(\pi-\theta) \\
    &\leq \pi \left((\pi-\theta)^{n-1} - \theta^{n-1} -\left(1-\frac{2\theta}{\pi}\right) \pi^{n-1}\right).
\end{align*}
As the right-hand side in the above expression is zero for $n=2$, an induction argument implies that the right-hand side is non-positive for any $2\leq n \in \N$. So from~\eqref{eq:derivative of f}, we see that $f_\theta^{(n)}(0) \leq 0$ for any $n \in \N$. Consequently, we can complete the proof by showing that for some $n\in \N$, the derivative $f_\theta^{(n)}$ is non-positive. To this end, note that for $n$ large enough we have $(\pi-\theta)^n/\pi^n \leq 1-2\theta/\pi$ (again because $\theta<\pi/2$). Thus 
\begin{align*}
    f^{(2n)}(x) &= (\pi-\theta)^{2n} \cosh((\pi-\theta)x) - (1-\theta/2\pi) \pi^{2n} \cosh(\pi x)- \theta^{2n} \cosh(\theta x)\\
    &\leq (\pi-\theta)^{2n} \cosh(\pi x) - \left(1-\frac{2\theta}{\pi}\right) \cosh(\pi x) \leq 0  \quad \mbox{for any $x\in \R$,} 
\end{align*}
which completes the proof of item~\ref{it:combinedsupremumnorm}.

Item~\ref{it:reproducingkernel} follows from Lemma~\ref{lem:Hardy} and the Fourier transform formula (see~\eqref{eq:exponentialfourier}) as
\begin{align*}
    \frac{2\pi}{\sin(\sigma)} \inner{\mathscr{T}_{z_0}\widehat{\mathscr{P}}_\sigma, h}_\sigma &= \frac{2\pi}{\sin(\sigma)} \int_\R \overline{\mathscr{P}_\sigma(x)\ee^{\ii x \overline{z_0}}} \left(\cosh(x)-\cos(\sigma)\right) \widecheck{h}(x)  \mathrm{d} x \\
    &= \int_\R \widecheck{h}(x)\ee^{-\ii xz_0} \mathrm{d}x = h(z_0).
\end{align*}
\end{proof}

We now turn to the proof of the residue theorem in Lemma~\ref{lem:residuethm}.
\begin{proof}[Proof of Lemma~\ref{lem:residuethm}]
Since $F$ has a simple pole at the origin, the function 
\begin{align}
    F(z) - \mathrm{Res}(F,0)/z \label{eq:poleexpansion}
\end{align}
is uniformly bounded in $\mathcal{S}^\ast$. As the residue is purely imaginary, the real part along the real axis can be written as
\begin{align*}
    \mathrm{Re} \left(F(x)\right) = \mathrm{Re} \left(F(x) - \frac{\mathrm{Res}(F,0)}{x}\right),
\end{align*}
and is therefore uniformly bounded in $x$.

To prove~\eqref{eq:residuethm}, we first claim that it suffices to consider $h= \widehat{\phi}$ for $\phi \in C_c^\infty(\R)$. To prove this claim, first note that $C_c^\infty(\R)$ is dense in the dual weighted space $\mathcal{L}^{2^\ast}_\sigma$ because the map
\begin{align*}
    g \in \mL^2(\R) \mapsto \Phi_\sigma g = \frac{g(\tau)}{\sqrt{\cosh(\tau)-\cos(\sigma)}} \in \mathcal{L}^{2^\ast}_\sigma
\end{align*}
is an isometry and maps $C_c^\infty(\R)$ bijectively to itself. Therefore the set $\{ \widehat{\phi} : \phi \in C_c^\infty(\R)\} \subset \mathbb{H}^2(\mathcal{S}^\ast)$ is dense by Lemma~\ref{lem:Hardy}. As convergence in $\mathbb{H}^2(\mathcal{S}^\ast)$ implies $\mL^2(\R)$ convergence along horizontal lines, and point evaluation is continuous in $\mathbb{H}^2(\mathcal{S}^\ast)$ (which can be seen from the reproducing kernel property in Lemma~\ref{lem:Poissonproperties}~\ref{it:reproducingkernel}), the claim now follows by approximating any $h\in \mathbb{H}^2(\mathcal{S}^\ast)$ by $\widehat{\phi}$ with $\phi \in C_c^\infty(\R)$ and passing to the limit in~\eqref{eq:residuethm}. Note that passing to the limit is justified by dominated convergence as $F$ is uniformly bounded away from zero and $\mathrm{Re}\, F(x)$ is bounded in $\R$.

To prove~\eqref{eq:residuethm} for functions of the form $h=\widehat{\phi}$, we apply Cauchy's integral theorem to the holomorphic function
\begin{align*}
    f(z) = F(z) h(z) \overline{h(\overline{z})}
\end{align*}
along the contour $\mathcal{C}_{\delta,R}$ composed of the parts 
\begin{align*}
    U_{\delta,R} &\coloneqq \left\{ x+\frac{1-\delta}{2}\ii: x \in \R, |x| \leq R\right\}, \quad L_R \coloneqq \{\pm R + \ii y: 0\leq y \leq 1-\delta\} \\
    D_{R,\delta} &\coloneqq \{ x\in \R: \delta < |x|\leq R\} \quad \hspace{2mm}\mbox{and} \hspace{2mm}\quad S_{\delta} \coloneqq \{ \delta \ee^{\ii \theta} : \theta \in (0,\pi)\},
\end{align*}
where $R>0$ and $0<\delta <1$. Then, since $\widehat{\phi}$ decays fast as $|x| \ra \infty$ and uniformly in $|y| \leq 1/2$ and $F$ is bounded, we can take the limit $R \rightarrow \infty$ to obtain
\begin{align*}
    \mathrm{Re} \left\{\int_\R F\left(x+\frac{1-\delta}{2} \ii\right) h_{\frac{1-\delta}{2}}(x) \overline{h_{-\frac{1-\delta}{2}}(x)}  \mathrm{d} x\right\}  &= \mathrm{Re} \left\{ \int_{|x| \geq \delta} F(x) |h(x)|^2 \mathrm{d} x - \int_{S_\delta} F(z) h(z) \overline{h(\overline{z})} \mathrm{d} z\right\}\\
    &=  \int_{|x| \geq \delta} \mathrm{Re} \{ F(x)\} |h(x)|^2 \mathrm{d} x - \mathrm{Re} \left\{ \int_{S_\delta} F(z) h(z) \overline{h(\overline{z})} \mathrm{d} z \right\}.
\end{align*}
To conclude, we can now take the limit $\delta \downarrow 0^+$ and use the expansion in~\eqref{eq:poleexpansion} to obtain~\eqref{eq:residuethm}.

\end{proof}

Finally, we prove the localization Lemma~\ref{lem:localizationlemma}.

\begin{proof}[Proof of Lemma~\ref{lem:localizationlemma}]
First, we note that the function $(mf)(x+\ii y)$ defined in~\eqref{eq:moperator} is simply the convolution against the kernel
\begin{align*}
    k_y(x) = 2\widehat{\mathscr{P}}_{\pi/2}(2x+2\ii y).
\end{align*}
In particular, we have
\begin{align}
    \mathcal{F}^{-1}(mf)(\tau) = 2\pi \widecheck{k}_0(\tau) \widecheck{f}(\tau) = 2\pi \mathscr{P}_{\pi/2}(\tau/2) \widecheck{f}(\tau) =  \frac{\widecheck{f}(\tau)}{\cosh(\tau/2)} \in \mathcal{L}^{2^\ast}_{\pi/2}(\R). \label{eq:importantid}
\end{align}
Hence, $mf \in \mathbb{H}^2(\mathcal{S}^\ast)$ by Lemma~\ref{lem:Hardy}. Moreover, as 
\begin{align*}
    \frac{\widecheck{h}_{1/2}(\tau) + \widecheck{h}_{-1/2}(\tau)}{2} = \cosh(\tau/2) \widecheck{h}(\tau),
\end{align*}
equation~\eqref{eq:firsttransformidentity} follows from~\eqref{eq:importantid}.

To prove the first part of~\eqref{eq:secondtransformidentity}, we note that, by~\eqref{eq:importantid} and the dilation property $\mathcal{F}^{-1} D_\delta f = D_{1/\delta} \mathcal{F}^{-1} f$ of the Fourier transform $\mathcal{F}$,
\begin{align*}
    \mathcal{F}^{-1} (m D_\delta f)(\tau) = \frac{\widecheck{f}(\delta \tau) \delta^{\frac12}}{\cosh(\tau/2)}.
\end{align*}
Hence, by Plancherel, re-scaling, and dominated convergence,
\begin{align*}
    \norm{(mD_\delta f)_0}_{\mL^2(\R)}^2 = 2\pi \int_\R \frac{|\widecheck{f}(\delta \tau)|^2 \delta}{\cosh(\tau/2)^2} \mathrm{d} \tau = 2\pi \int_\R \frac{|\widehat{f}(\tau)|^2}{\cosh\left(\tau/(2\delta)\right)^2} \mathrm{d} \tau \rightarrow 0 \quad \mbox{as $\delta \downarrow 0^+$.}
\end{align*}
To prove the second statement in~\eqref{eq:secondtransformidentity}, we note that
\begin{align*}
\frac{\mathcal{F}^{-1} D_{1/\delta} (m D_\delta f)_{1/2}(\tau) - \mathcal{F}^{-1} D_{1/\delta} (m D_\delta f)_{-1/2}(\tau)}{2} = \frac{\sinh(\tau/(2\delta))}{\cosh(\tau/(2\delta))} \widecheck{f}(\tau).
\end{align*}
As $\tanh(\tau/(2\delta)) \ra \mathrm{sign}(\tau)$ for every $\tau \neq 0$, the result follows from Plancherel's theorem, dominated convergence, and the well-known identity for the Hilbert transform (cf. \cite[Equation 5.3]{Kin09})
\begin{align*}
    \mathcal{F}^{-1}(Hf)(\tau) = \ii \mathrm{sign}(\tau) \widecheck{f}(\tau).
\end{align*}
\end{proof}

\section*{Disclosure statement}
The authors report there are no competing interests to declare.

\bigskip
\end{document}